\documentclass{amsart}
\usepackage{amsfonts,amssymb}
\usepackage{epsfig}
\newtheorem{theorem}{Theorem}
\newtheorem{lemma}{Lemma}
\newtheorem{corollary}{Corollary}

\newcommand{\integers}{{\mathbb Z}}
\newcommand{\realnos}{{\mathbb R}}

\def\diag{{\rm diag}}
\newcommand{\p}{\phantom}

\begin{document}

\title[Integral Congruence Two Hyperbolic 5-Manifolds]
{Integral Congruence Two Hyperbolic 5-Manifolds}

\author{John G. Ratcliffe and Steven T. Tschantz}

\address{Department of Mathematics, Vanderbilt University, Nashville, TN 37240
\vspace{.1in}}

\email{ratclifj@math.vanderbilt.edu}


\date{}


\begin{abstract}
In this paper, we classify all the orientable hyperbolic 5-manifolds 
that arise as a hyperbolic space form $H^5/\Gamma$ 
where $\Gamma$ is a torsion-free subgroup of minimal index 
of the congruence two subgroup $\Gamma^5_2$   
of the group $\Gamma^5$ of positive units of the Lorentzian quadratic form 
$x_1^2+\cdots+x_5^2-x_6^2$. 
We also show that $\Gamma^5_2$ is a reflection group 
with respect to a 5-dimensional right-angled convex polytope in $H^5$. 
As an application, we construct a hyperbolic 5-manifold 
of smallest known volume $7\zeta(3)/4$.
\end{abstract}

\maketitle

\section{Introduction} 
According to Carl Ludwig Siegel,
\begin{quote}
\em the mathematical universe is inhabited not only by important 
species but also by interesting individuals. 
\end{quote} 
Arithmetic hyperbolic manifolds form an important species of hyperbolic manifolds. 
In this paper, we construct and classify some interesting individuals. 
These individuals are interesting because of their dimension, five, 
and the value, $28\zeta(3)$, of their volume. 
The volumes of our 5-manifolds are interesting 
because the value of the Riemann zeta function at three, 
$\zeta(3)$, is known to be an irrational number [11], 
in contrast to dimension three, where the irrationality 
of volumes of hyperbolic 3-manifolds is unknown.

Our manifolds are explicitly described by gluing together the sides 
of a hyperbolic convex polytope. 
These manifolds are the first examples of explicitly described 
complete hyperbolic 5-manifolds of finite volume. 
One of our examples has a group of symmetries of order 16 
that acts freely on the manifold. 
The resulting quotient hyperbolic 5-manifold has 
volume, $7\zeta(3)/4$, which is the smallest known volume of 
a complete hyperbolic 5-manifold. 
The previous smallest known volume was $14\zeta(3)$ 
which is the volume of a hyperbolic 5-manifold constructed 
by B. Everitt [1]. 
A volume lower bound of R. Kellerhals [3] suggests 
that our quotient manifold has very small volume 
and may be a minimum volume open hyperbolic 5-manifold. 
See our paper [9], for the volumes 
of some arithmetic hyperbolic $5$-manifolds. 

We now set up notation. 
A real $(n+1)\times (n+1)$ matrix $A$ is said to be {\it Lorentzian}
if $A$ preserves the {\it Lorentzian inner product}
$$x \circ y\, =\, x_1y_1 + x_2y_2 + \cdots + x_ny_n - x_{n+1}y_{n+1}.$$
The {\it hyperboloid model} of hyperbolic $n$-space is the metric space
$$H^n\, =\, \{x \in \realnos^{n+1}: \hbox{$x \circ x = -1$ and $x_{n+1} >
0$}\}$$
with metric $d$ defined by
$$\cosh d(x,y) = -x \circ y.$$
A Lorentzian $(n+1)\times (n+1)$ matrix $A$ is said to be
either {\it positive} or {\it negative}
according as $A$ maps $H^n$ to $H^n$ or $-H^n$.
The isometries of $H^n$ correspond to the positive Lorentzian
$(n+1)\times (n+1)$ matrices.

Let $\Gamma^n$ be the group of positive
Lorentzian $(n+1) \times (n+1)$ matrices with integer entries.
The group $\Gamma^n$ is an infinite discrete subgroup
of the group ${\rm O}(n,1)$ of Lorentzian  $(n+1) \times (n+1)$ matrices.
The {\it principal congruence two subgroup} of $\Gamma^n$
is the group $\Gamma^n_2$ of all matrices in $\Gamma^n$ that are congruent
to the identity matrix modulo two.
The congruence two subgroup $\Gamma^n_2$ is not torsion-free, but
it only has 2-torsion by Theorem IX.7 of Newman [5]. 

\section{Congruence Two Subgroup} 

In this section, we determine the structure of the congruence two
subgroup $\Gamma^5_2$ of the group $\Gamma^5$ of
integral, positive, Lorentzian $6 \times 6$ matrices. 
Before reading this section, we advise the reader to read the first four sections 
of our paper [10] which determines the structure of 
$\Gamma^2_2, \Gamma^3_2$, and $\Gamma^4_2$.  
Vinberg [12] proved that the group $\Gamma^5$
is a reflection group with respect to a noncompact
5-simplex $\Delta^5$ in $H^5$ whose Coxeter diagram is
\begin{center}
\begin{picture}(210,80)(-80,-50)
\put(-48,7){3}
\put(-75,0){\circle*{4}}
\put(-75,0){\line(1,0){50}}
\put(-2,7){3}
\put(-25,0){\line(1,0){50}}
\put(-25,0){\circle*{4}}
\put(25,0){\circle*{4}}
\put(48,7){3}
\put(25,0){\line(1,0){50}}
\put(75,0){\circle*{4}}
\put(32,-26){3}
\put(25,0){\line(0,-1){48}}
\put(25,-48){\circle*{4}}
\put(75,0){\line(1,0){50}}
\put(75,0){\circle*{4}}
\put(125,0){\circle*{4}}
\put(98,7){4}
\end{picture}
\end{center}
The volume of $\Delta^5$ is $7\zeta(3)/15360$, see [9].  
Vertices for $\Delta^5$ are
$$\begin{array}{l}
(0,0,0,0,0,1),\,
(\sqrt{6}/6,\!\sqrt{6}/6,\!\sqrt{6}/6,0,0,\sqrt{6}/2),\\
(\sqrt{2}/2,\!\sqrt{2}/2,0,0,0,\sqrt{2}),\,
(\sqrt{5}/5,\sqrt{5}/5,\sqrt{5}/5,\sqrt{5}/5,0,3\sqrt{5}/5),\\
(1/2, 1/2, 1/2, 1/2, 1/2, 3/2),\,
(1,0,0,0,0,1)\ \hbox{(at infinity)}.
\end{array}$$
The group $\Gamma^5$ is generated by the following six matrices 
that represent the reflections of $H^5$ in the sides of $\Delta^5$, 
$$\left(\begin{array}{cccccc}
          0 & 1 & 0 & 0 & 0 & 0 \\
          1 & 0 & 0 & 0 & 0 & 0 \\
          0 & 0 & 1 & 0 & 0 & 0 \\
		  0 & 0 & 0 & 1 & 0 & 0 \\
		  0 & 0 & 0 & 0 & 1 & 0 \\
		  0 & 0 & 0 & 0 & 0 & 1 
        \end{array} \right),
\left(\begin{array}{cccccc}
          1 & 0 & 0 & 0 & 0 & 0 \\
          0 & 0 & 1 & 0 & 0 & 0 \\
          0 & 1 & 0 & 0 & 0 & 0 \\
		  0 & 0 & 0 & 1 & 0 & 0 \\
		  0 & 0 & 0 & 0 & 1 & 0 \\
		  0 & 0 & 0 & 0 & 0 & 1 
        \end{array} \right),
\left(\begin{array}{cccccc}
          1 & 0 & 0 & 0 & 0 & 0 \\
          0 & 1 & 0 & 0 & 0 & 0 \\
          0 & 0 & 0 & 1 & 0 & 0 \\
		  0 & 0 & 1 & 0 & 0 & 0 \\
		  0 & 0 & 0 & 0 & 1 & 0 \\
		  0 & 0 & 0 & 0 & 0 & 1 
        \end{array} \right),$$
$$\left(\begin{array}{cccccc}
          \!\!\!\phantom{-}0 &           \!-1 &           \!-1 & 0 & 0 & 1 \\
                   \!\! \!-1 & \!\phantom{-}0 &           \!-1 & 0 & 0 & 1 \\
					\!\!\!-1 &           \!-1 & \!\phantom{-}0 & 0 & 0 & 1 \\
         \! \!\!\phantom{-}0 & \!\phantom{-}0 & \!\phantom{-}0 & 1 & 0 & 0 \\
         \! \!\!\phantom{-}0 & \!\phantom{-}0 & \!\phantom{-}0 & 0 & 1 & 0 \\          
			       \! \!\!-1 &           \!-1 &           \!-1 & 0 & 0 & 2
       \end{array} \right),
\left(\begin{array}{cccccc}
          1 & 0 & 0 & 0 & 0 & 0 \\
          0 & 1 & 0 & 0 & 0 & 0 \\
          0 & 0 & 1 & 0 & 0 & 0 \\
		  0 & 0 & 0 & 0 & 1 & 0 \\
		  0 & 0 & 0 & 1 & 0 & 0 \\
		  0 & 0 & 0 & 0 & 0 & 1 
        \end{array} \right), 
\left(\begin{array}{cccccc}
          1 & 0 & 0 & 0 & \!\phantom{-}0 & 0 \\
		  0 & 1 & 0 & 0 & \!\phantom{-}0 & 0 \\
		  0 & 0 & 1 & 0 & \!\phantom{-}0 & 0 \\
		  0 & 0 & 0 & 1 & \!\phantom{-}0 & 0 \\            
		  0 & 0 & 0 & 0 &		    \!-1 & 0 \\
		  0 & 0 & 0 & 0 & \!\phantom{-}0 & 1
       \end{array} \right).$$

Let $\Sigma^5$ be the group generated the first five matrices in
the above list of matrices.
Then $\Sigma^5$ is the group generated
by the reflections in the sides of $\Delta^5$ incident with
the vertex $(1/2,1/2,1/2,1/2,1/2,3/2)$. 
Therefore $\Sigma^5$ is isomorphic to a spherical 4-simplex reflection group
whose Coxeter diagram is obtained from the Coxeter diagram of
$\Gamma^5_2$ by deleting its right most vertex
and its adjoining edge. Hence $\Sigma^5$ has the Coxeter diagram
\begin{center}
\begin{picture}(160,80)(-80,-50)
\put(-48,7){3}
\put(-75,0){\circle*{4}}
\put(-75,0){\line(1,0){50}}
\put(-2,7){3}
\put(-25,0){\line(1,0){50}}
\put(-25,0){\circle*{4}}
\put(25,0){\circle*{4}}
\put(48,7){3}
\put(25,0){\line(1,0){50}}
\put(75,0){\circle*{4}}
\put(32,-26){3}
\put(25,0){\line(0,-1){48}}
\put(25,-48){\circle*{4}}
\end{picture}
\end{center}
The order of the Coxeter group $\Sigma^5$ is 1920.
None of the elements of $\Sigma^5$ are in $\Gamma^5_2$. 
By Lemma 16 of [9], 
the index of $\Gamma_2^5$ in $\Gamma^5$ is 1920.  
Hence $\Sigma^5$ is a set of coset representatives
for $\Gamma_2^5$ in $\Gamma^5$.  We therefore
have a natural, split, short, exact sequence of groups
$$ 1 \to \Gamma_2^5 \to \Gamma^5 \to \Sigma^5 \to 1.$$

Let $P^5 = \Sigma^5\Delta^5$.  
Then $P^5$ is a convex hyperbolic polytope  
which is subdivided into 1920 copies of $\Delta^5$ 
obtained by reflecting in the sides of $\Delta^5$ 
that are incident with the vertex $(1/2,1/2,1/2,1/2,1/2,3/2)$. 
The polytope $P^5$ has 16 actual vertices and 10 ideal vertices. 
The ideal vertices are
$$\begin{array}{l} 
 (1, 0, 0, 0, 0, 1), 
 (0, 1, 0, 0, 0, 1), 
 (0, 0, 1, 0, 0, 1), 
 (0, 0, 0, 1, 0, 1),
 (0, 0, 0, 0, 1, 1),\\
 (0, 1, 1, 1, 1, 2), 
 (1, 0, 1, 1, 1, 2), 
 (1, 1, 0, 1, 1, 2), 
 (1, 1, 1, 0, 1, 2), 
 (1, 1, 1, 1, 0, 2).
 \end{array}$$
The 10 ideal vertices of $P^5$ are the vertices of a regular ideal 
5-dimensional cross polytope $R^5$ centered at the point 
$(1/2,1/2,1/2,1/2,1/2,3/2)$. 
The regular cross polytope $R^5$ has 32 sides 
each of which is a regular ideal 4-simplex.  

Let $e_1,\ldots, e_6$ be the standard basis of $\realnos^6$ and let 
$T^5$ be the 5-simplex in $H^5$ with one actual vertex at 
the origin $e_6$ of $H^5$ and five ideal vertices 
$e_i+e_6$ for $i=1,\ldots,5$. 
We call $T^5$ a {\it corner} 5-{\it simplex}, with $e_6$
the {\it corner vertex} of $T^5$, and each side of $T^5$ incident
with $e_6$ a {\it corner side} of $T^5$. 
The 5-simplex $T^5$ has five corner sides. 
The sixth side of $T^5$ is a regular ideal 4-simplex. 
The dihedral angle between any two corner sides of $T^5$ is a right angle. 
The dihedral angle between the ideal side of $T^5$ and a corner side is $\pi/3$. 
Thus $T^5$ is a hyperbolic Coxeter simplex. 

The regular ideal cross polytope $R^5$ is subdivided into  
32 copies of $T^5$ fitting together around their corner 
vertices at the point $(1/2,1/2,1/2,1/2,1/2,3/2)$. 
Each side of $R^5$ is the ideal side of one 
of the 32 corner simplices subdividing $R^5$. 
The dihedral angles between the ideal sides and a common corner side 
of two adjacent corner simplices subdividing $R^5$ add up 
to the dihedral angle $2\pi/3$ of $R^5$. 

The polytope $P^5$ is obtained from $R^5$ by gluing onto 
each of 16 mutually nonadjacent sides of $R^5$ a corner 5-simplex   
with one of the attached corner 5-simplices being $T^5$. 
The corner vertices of the attached corner 5-simplices  
are the actual vertices of $P^5$.  
The actual vertices of $P^5$ are 
$$\begin{array}{l}
(0, 0, 0, 0, 0, 1), (0, 0, 1, 1, 1, 2), (0, 1, 0, 1, 1, 2), (0, 1, 1, 0, 1, 2), \\
(0, 1, 1, 1, 0, 2), (1, 0, 0, 1, 1, 2), (1, 0, 1, 0, 1, 2), (1, 0, 1, 1, 0, 2), \\ 
(1, 1, 0, 0, 1, 2), (1, 1, 0, 1, 0, 2), (1, 1, 1, 0, 0, 2), (1, 1, 1, 1, 2, 3), \\
(1, 1, 1, 2, 1, 3), (1, 1, 2, 1, 1, 3), (1, 2, 1, 1, 1, 3), (2, 1, 1, 1, 1, 3). 
\end{array}$$ 

The polytope $P^5$ has 16 sides. 
The remaining 16 sides of $R^5$ extend to the 16 sides of $P^5$. 
Each of these 16 sides of $R^5$ is a regular ideal 4-simplex  
which is coplanar with a corner side of five of the attached  
corner 5-simplices. Thus each side of $P^5$ is obtained 
from a regular ideal 4-simplex $S^4$ by 
gluing on a corner 4-simplex onto each side of $S^4$. 
Hence, each side of $P^5$ is congruent to the polytope $P^4$ 
described in our paper [10]. 
The group of symmetries of $P^5 = \Sigma^5\Delta^5$ is $\Sigma^5$ 
which is a subgroup of index two of the group of symmetries of $R^5$.

The dihedral angles of $P^5$ are the dihedral angles 
of the attached corner 5-simplices between corner sides, 
and so $P^5$ is a right-angled polytope. 
It is worth noting that L. Potyagailo and E. Vinberg [6] 
have proven that a 5-dimensional, right-angled, hyperbolic polytope  
of finite volume has at least 16 sides, 
and so $P^5$ is a 5-dimensional, 
right-angled, hyperbolic polytope of finite volume 
with the smallest possible number of sides. 
 
The polytope $P^5$ has 16 sides, five of them are coplaner 
with the coordinate hyperplanes $x_i=0$, for $i=1,\ldots, 5$. 
Let $D = {\rm diag}(1,1,1,1,-1,1)$ be the last generator of $\Delta^5$. 
The diagonal matrix $D$ represents the reflection of $H^5$ 
in the hyperplane $x_5=0$. The matrix $D$ is in $\Gamma^5_2$. 
The reflections of $H^5$ in the 16 sides of $P^5$ are represented 
by matrices of the form $ADA^{-1}$ with $A$ in $\Sigma^5$. 
Each such matrix $ADA^{-1}$ is in $\Gamma^5_2$, 
since $\Gamma^5_2$ is a normal subgroup of $\Gamma^5$. 
Now since $\Sigma^5$ is a set of coset representatives
for $\Gamma_2^5$ in $\Gamma^5$, we have that $P^5 = \Sigma^5\Delta^5$
is a fundamental polytope for $\Gamma_2^5$.
We therefore have the following theorem.

\begin{theorem} 
The congruence two subgroup $\Gamma_2^5$ of the group
$\Gamma^5$ of integral, positive, Lorentzian $6 \times 6$ matrices is a
reflection group with respect to the right-angled polytope $P^5$.
\end{theorem}

We now describe the 16 matrices that represent the reflections of $H^5$ 
in the sides of $P^5$. 
Five of the sides of $P^5$ are coplanar 
with the coordinate hyperplanes $x_i=0$, for $i=1,\ldots, 5$. 
The reflections in the five coordinate hyperplanes 
are represented by diagonal matrices with diagonal entries all 1 except for 
a $-1$ in the $i$th diagonal slot. 
The polytope $P^5$ has 10 sides each of which is perpendicular 
to three of the coordinate hyperplanes in all 10 combinations of 
three coordinate hyperplanes. 
The matrices representing the reflections in these 10 sides are obtained 
from the matrix 
$$\left(\begin{array}{cccccc}
          1 & 0 & 0 & \!\phantom{-}0 & \!\phantom{-}0 & 0 \\
          0 & 1 & 0 & \!\phantom{-}0 & \!\phantom{-}0 & 0 \\
          0 & 0 & 1 & \!\phantom{-}0 & \!\phantom{-}0 & 0 \\
		  0 & 0 & 0 &			\!-1 &           \!-2 & 2 \\
          0 & 0 & 0 &           \!-2 &           \!-1 & 2 \\
          0 & 0 & 0 &           \!-2 &           \!-2 & 3 \\
       \end{array} \right)$$
by conjugating by a permutation matrix of the first five coordinates. 
Note that the three diagonal ones can lie in any one of the 
10 combinations of three diagonal positions in the first five 
diagonal slots. 
The remaining side of $P^5$ is not incident with any of the coordinate 
hyperplanes. We call this side the {\it far side} of $P^5$.  
The reflection of $H^5$ in the far side of $P^5$ is represented by 
the matrix 
$$\left(\begin{array}{cccccc}
          -1 & -2 & -2 & -2 & -2 & 4 \\
          -2 & -1 & -2 & -2 & -2 & 4 \\
          -2 & -2 & -1 & -2 & -2 & 4 \\
		  -2 & -2 & -2 & -1 & -2 & 4 \\
          -2 & -2 & -2 & -2 & -1 & 4 \\
          -4 & -4 & -4 & -4 & -4 & 9 \\
       \end{array} \right).$$

Let $K^5$ be the group of 32 diagonal matrices  
${\rm diag}(\pm 1,\pm 1,\pm ,\pm 1,\pm 1, 1).$ 
Then $K^5$ is a subgroup of $\Gamma^5_2$. 
The next corollary follows immediately from Theorem 1. 
\begin{corollary} 
Every nonidentity element of $\Gamma_2^5$ of finite order has order two,
every finite subgroup of $\Gamma_2^5$ is conjugate in $\Gamma^5$
to a subgroup of the elementary 2-group $K^5$, and
there are 16 conjugacy classes of maximal finite subgroups of $\Gamma_2^5$
in $\Gamma_2^5$ corresponding to the 16 actual vertices of $P^5$.
\end{corollary}

\section{Torsion-Free Subgroups of Index 32} 

Let $\Gamma$ be a torsion-free subgroup of $\Gamma^5_2$ of finite index.
Then the finite group $K^5$ acts freely on the set of cosets of $\Gamma$
in $\Gamma^5_2$ by $\Gamma g \mapsto \Gamma gk$,
since $\Gamma$ is torsion-free.
Therefore $|K^5| = 32$ divides $[\Gamma^5_2:\Gamma]$.

Now suppose that $[\Gamma^5_2:\Gamma] = 32$.  
Then the convex polytope
$Q^5 = K^5P^5$ is a fundamental polytope for $\Gamma$. 
The polytope is subdivided into 32 copies of $P^5$, one of which is $P^5$,  
that fit together around the origin $e_6$ of $H^5$, 
with one copy of $P^5$ in each of the 32 orthants of $H^5$ bounded by 
five coordinate hyperplanes of $H^5$. 
 
The polytope $Q^5$ has 72 sides of two different types. 
A side of $P^5$ lying in a coordinate hyperplane is identified 
with a side of the copy of $P^5$ across the coordinate hyperplane. 
This eliminates $5\cdot 32$ sides of the 32 copies of $P^5$ subdividing $Q^5$. 
Let $S$ be a side of $P^5$ that is perpendicular to three coordinate hyperplanes 
and let $K_S$ be the subgroup of $K^5$ generated by the reflections  
in the three coordinate hyperplanes perpendicular to $S$.  
Then $K_S$ leaves the hyperplane of $H^5$ spanned by $S$ invariant. 
Hence the eight sides $\{kS: k \in K_S\}$ of the copies of $P^5$ 
subdividing $Q^5$ are coplanar. 
These eight copies of $P^4$ form a side $K_SS$ of $Q^5$. 
The intersection of the eight copies of $P^4$ forming $K_SS$ 
is a hyperbolic line which we call the {\it axis} of the side $K_SS$ of $Q^5$. 
The axis of the side $K_SS$ is an edge of $S$ joining two ideal vertices 
$e_i+e_6$ and $e_j+e_6$ of $S$. 
There are a total of 40 sides of $Q^5$ of this type 
which accounts for $40\cdot 8 = 10\cdot 32$ sides 
of the 32 copies of $P^5$ subdividing $Q^5$. 
Now $K^5$ maps the far side of $P^5$ to 32 sides of $Q^5$ 
that are congruent to $P^4$. 
This accounts for all of the $16\cdot 32$ sides 
of the 32 copies of $P^5$ subdividing $Q^5$. 
Thus $Q^5$ has a total of 72 sides, 32 small sides congruent to $P^4$, 
and 40 large sides made up of eight copies of $P^4$ glued together 
around a common infinite edge. 
The axes of the 40 large sides of $Q^5$ are the edges 
of the regular ideal 5-dimensional cross polytope $O^5$ 
with vertices $\pm e_i + e_6$, for $i=1,\ldots, 5$. 
There is one small side of $Q^5$ in each of the 32 orthants of $H^5$ 
bounded by five coordinate hyperplanes. 
The small sides of $Q^5$ are pairwise disjoint. 
The polytope $Q^5$ has the same symmetry group 
as the regular cross polytope $O^5$. 
All the dihedral angles of the sides of $Q^5$ are right angles and 
all the dihedral angles of $Q^5$ are right angles.   

Let $S$ be a side of $Q^5$ and let $r$ be the reflection 
of $H^n$ in the side $S$. 
Then there is a reflection $s$ in a side of $P^5$ and an element $\ell$ 
of $K^5$ such that $r=\ell s \ell^{-1}$. 
Hence $r$ is in $\Gamma^5_2$. 
Now since $K^5$ is a set of coset representatives for $\Gamma$ in $\Gamma^5_2$, 
there is an element $g$ of $\Gamma$ 
and an element $k$ of $K^5$ such that $r=gk$.   
Then $g=rk$. Thus $rk$ is a side-pairing 
transformation in $\Gamma$ for $Q^5$ 
that pairs the side $S'=kS$ of $Q^5$ to the side $S$ of $Q^5$. 
Summarizing, we have the following lemma. 

\begin{lemma} 
Let $\Gamma$ be a torsion-free subgroup of $\Gamma^5_2$ of index 32. 
Then the right-angled polytope $Q^5$ is an exact fundamental 
polytope for $\Gamma$. 
The side-pairing transformations of $Q^5$ 
are of the form $rk$ where $k$ is in $K^5$ 
and $r$ is the reflection in a side of $Q^5$.
\end{lemma}

\section{Geometry of the Fundamental Polytope $Q^5$} 

The polytope $Q^5$ has 160 actual vertices and 90 ideal vertices. 
The polytope $P^4$ has five actual vertices and five ideal vertices. 
All the actual vertices of $Q^5$ are vertices of small sides 
for a total of $5\cdot 32 = 160$ actual vertices. 
The 160 actual vertices are obtained from $(\pm 2,\pm 1,\pm 1,\pm 1,\pm 1, 3)$ 
by permuting the first five coordinates. 
The ten ideal vertices $\pm e_i+e_6$ of the regular cross polytope 
$O^5$ are ideal vertices of $Q^5$, which we call {\it large ideal vertices}. 
Each large ideal vertex is a vertex of eight large sides and no small sides. 
All the other ideal vertices of $Q^5$ are ideal vertices of small sides, 
with each ideal vertex shared by two small sides in adjacent orthants, 
for a total of $(5\cdot 32)/2 = 80$ ideal vertices of small sides, 
which we call {\it small ideal vertices}. 
The 80 small ideal vertices are obtained from $(0,\pm 1,\pm 1,\pm 1,\pm 1, 2)$ 
by permuting the first five coordinates. 

We next describe the ridges (3-faces) of $Q^5$. 
The polytope $Q^5$ has 560 ridges of two different types. 
The polytope $P^4$ has 10 sides, each of which 
is congruent to the polyhedron $P^3$ described in our paper [10]. 
See Figure 1.  
The polyhedron $P^3$ has five vertices, two actual and three ideal, 
and six sides, each of which is an ideal right triangle congruent to $P^2$. 
All the dihedral angles of $P^3$ are right angles. 
\begin{figure}[b] 
\epsfxsize=70mm
\epsfbox{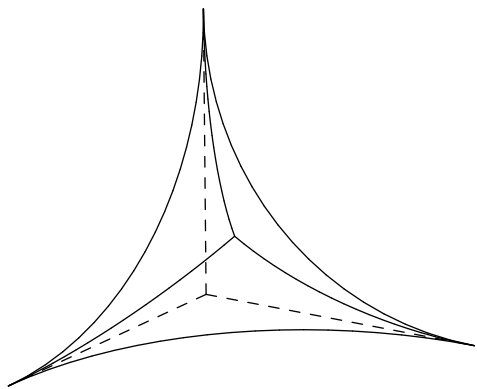}
\caption{The polyhedron $P^3$}
\end{figure}
The polytope $Q^5$ has $10\cdot 32=320$ small ridges congruent to $P^3$. 
Each small ridge is a side of just one small side of $Q^5$ 
and just one large side of $Q^5$.  

Let $S$ and $T$ be adjacent sides of $P^5$ that are both perpendicular 
to three coordinate hyperplanes. 
These two sets of three coordinate hyperplanes have two hyperplanes in common. 
The intersection of $S$ and $T$ is a ridge 
$R$ of $P^5$, and so $R$ is congruent to $P^3$. 
Now $R$ is perpendicular to the two common coordinate hyperplanes and 
lies in the other two coordinate hyperplanes. 
Let $K_R$ be the subgroup of $K^5$ generated by the reflections 
in the two coordinate hyperplanes perpendicular to $R$. 
Then $K_R$ leaves the 3-plane of $H^5$ spanned by $R$ invariant. 
Hence the four ridges $\{kR: k\in K_R\}$ of copies of $P^5$ subdividing $Q^5$ 
are coplanar.  These four copies of $P^3$ form a ridge $K_RR$ of $Q^5$. 
The intersection of the four copies of $P^3$ forming $K_RR$ is a hyperbolic ray 
which we call the {\it axis} of the ridge $K_RR$ of $Q^5$. 
The axis of the ridge $K_RR$ is an edge of $R$ joining a large ideal vertex 
to an actual vertex of $R$. 
The polyhedron $K_RR$ is congruent to the top half of 
the rhombic dodecahedron $Q^3$ in Figure 5 of our paper [10]. 
See Figure 2. 
\begin{figure} 
\epsfxsize=100mm
\epsfbox{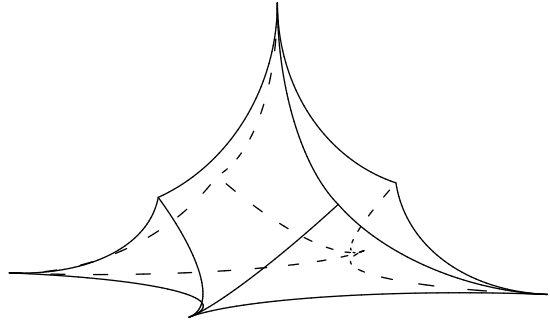}
\caption{A large ridge of $Q^5$}
\end{figure}
The polyhedron $K_RR$ has nine vertices, four actual and five ideal, 
and nine sides, four triangles congruent to $P^2$, 
four rhombi made up of two copies of $P^2$ 
glued together along a line edge, 
and an ideal square congruent to $Q^2$. 
All the dihedral angles of the polyhedron $K_RR$ are right angles. 
A large side of $Q^5$ has 20 sides, eight small sides in common with adjacent 
small sides of $Q^5$ and 12 large sides in common 
with adjacent large sides of $Q^5$. 
Thus $Q^5$ has a total of $40 \cdot 8 = 320$ small ridges congruent to $P^3$ 
and $(40 \cdot 12)/2 = 240$ large ridges made up of four copies of $P^3$ 
glued together around a common ray edge.  

Let $v$ be a large ideal vertex of $Q^5$ 
and let $T$ be a large side of $Q^5$ with ideal vertex $v$. 
Each large side of $T$, with ideal vertex $v$, is adjacent to four 
large sides of $T$, with ideal vertex $v$, 
along its four rhombic faces. 
Therefore the link of the ideal vertex $v$ in $T$ 
is a cube and $T$ has six large sides with ideal vertex $v$. 
The polytope $Q^5$ has eight large sides with ideal vertex $v$. 
Each large side of $Q^5$, with ideal vertex $v$, is adjacent 
to six large sides of $Q^5$ with ideal vertex $v$. 
Hence the link of the large ideal vertex $v$ of $Q^5$ is a 4-cube. 

Let $T$ be a large side of $Q^5$ with large ideal vertices $v$ and $w$. 
Then $T$ has six large sides with ideal vertex $v$   
and six large sides with ideal vertex $w$. 
Each large side of $T$, with ideal vertex $v$, 
is adjacent to a large side of $T$, with ideal vertex $w$, 
along a common ideal square face. 
Each large side of $T$ is adjacent to four small sides 
along its four triangle faces and each small side of $T$ 
is adjacent to six large sides of $T$ 
along its six triangle faces. 

Let $u$ be a small ideal vertex of $Q^5$ 
and let $S$ be a small side of $Q^5$ with ideal vertex $u$. 
The link of $u$ in $S$ is a cube. 
Hence the small side $S$ of $Q^5$ is adjacent 
to six large sides of $Q^5$ with ideal vertex $u$. 
Let $R$ be a large ridge of $Q^5$ with ideal vertex $u$. 
Then $u$ is a vertex of the ideal square face of $R$. 
The link of $u$ in $R$ is a rectangle made up 
of two adjacent squares. 
Let $T$ be a large side of $Q^5$ with ideal vertex $u$. 
Then two adjacent large sides of $T$, 
with a common ideal square face, 
have $u$ as an ideal vertex. 
Hence the link of $u$ in $T$ is a rectangular box 
made up of two adjacent cubes. 
Thus $T$ has four large sides and two small sides 
with ideal vertex $u$. 
Finally, the link of the small ideal $u$ in $Q^5$ 
is a rectangular 4-box made up of two adjacent 4-cubes. 
Thus $Q^5$ has six large sides and two small sides with ideal vertex $u$. 
The large side $T$ of $Q^5$ has $8\cdot 2 = 16$ actual vertices, 
two large ideal vertices, and $(8\cdot 3)/2 = 12$ small ideal vertices. 

The polytope $Q^5$ has 1360 2-faces of three different types, 
triangles congruent to $P^2$, rhombi made up of two copies of $P^2$, 
and ideal squares. 
Each triangle face of $Q^5$ is a face of a small side of $Q^5$. 
Now $P^3$ has six faces congruent to $P^2$ and $P^4$ has ten sides 
congruent to $P^3$.  Hence $P^4$ has 30 triangle faces. 
Therefore $Q^5$ has $32\cdot 30 = 960$ triangle faces. 
Let $F$ be a rhombic face of $Q^5$. 
Then $F$ is a face of a large ridge of $Q^5$,  
and so $F$ has a large ideal vertex $v$. 
Hence $F$ corresponds to an edge in the link of $v$. 
The link of $v$ is a 4-cube. 
A 4-cube has 32 edges. 
Therefore $Q^5$ has $10\cdot 32 = 320$ rhombic faces. 
Let $S$ be an ideal square face of $Q^5$. 
Then $S$ is a face of a large ridge $R_0$ of $Q^5$ 
with large ideal vertex $v$. 
Now $R_0$ is the intersection of two large sides 
$T_1$ and $T_2$ of $Q^5$ 
with large ideal vertices $u, v$ and $v, w$, respectively. 
The ridge $R_0$ is adjacent to a large ridge $R_i$ in $T_i$ 
along $S$ for $i= 1,2$. 
The ridges $R_1$ and $R_2$ are adjacent along $S$ 
in the large side $T_3$ with large ideal vertices $u$ and $w$. 
Thus $S$ is a face of exactly three large ridges in $Q^5$. 
Therefore there are a total of $240/3=80$ ideal square faces of $Q^5$. 

The polytope $Q^5$ has 1120 edges of two different types, rays and lines. 
Every line edge of $Q^5$ is an edge of a small side of $Q^5$. 
The line edges of $P^4$ are the edges of a regular ideal 4-simplex. 
Hence $P^4$ has 10 line edges. 
Therefore $Q^5$ has $32\cdot 10 = 320$ line edges. 
The link of an actual vertex $v$ of $Q^5$ is a right-angled 
spherical 4-simplex whose vertices correspond to the edges 
of $Q^5$ with vertex $v$. 
Thus $Q^5$ has $5\cdot 160 = 800$ ray edges. 

The Euclidean closure of $Q^5$ in the projective disk model of hyperbolic $5$-space 
is a Euclidean polytope $\overline{Q}\hbox{}^5$, 
and so the Euler characteristic of $\overline{Q}\hbox{}^5$ is one. 
As a consistency check, we compute $\chi(\overline{Q}\hbox{}^5)$ 
in terms of the face decomposition of $\overline{Q}\hbox{}^5$,  
$$\chi(\overline{Q}\hbox{}^5) = 250-1120+1360-560+72-1 = 1.$$

\section{Integral, Congruence Two, Orientable, Hyperbolic 5-Manifolds} 

In order to describe a torsion-free subgroup of $\Gamma^5_2$ of index 32, 
we have to construct a proper side-pairing of the polytope $Q^5$ 
with side-pairing transformation of the form in Lemma 1. 
See Sections 11.1 and 11.2 of Ratcliffe [7] for a discussion 
of complete gluing of hyperbolic $n$-manifolds. 
The normalized solid angle subtended by a point in the interior 
of a $k$-face of $Q^5$ is $2^{k-5}$ for each $k=0,\ldots,5$, 
since $Q^5$ is right-angled. 
Hence, in a proper side-pairing of $Q^5$, we must have $2^{5-k}$ points 
in each cycle of points in the interior of $k$-faces. 
Thus in a proper side-pairing of $Q^5$, 
there are $2^{5-k}$ $k$-faces of a single type 
in each cycle of $k$-faces for each $k=0,\ldots,5$. 
Hence in a proper side-pairing of $Q^5$ 
there are $160/32 = 5$ cycles of actual vertices,     
$800/16 = 50$ cycles of ray edges, 
$320/16 = 20$ cycles of line edges,  
$960/8=120$ cycles of triangle faces,  
$320/8=40$ cycles of rhombic faces,  
$80/8=10$ cycles of ideal squares, 
$320/4=80$ cycles of small ridges, 
$240/4=60$ cycles of large ridges, 
$32/2 = 16$ cycles of small sides, 
and $40/2=20$ cycles of large sides. 

Let $M$ be a complete hyperbolic 5-manifold 
obtained by gluing together pairs of sides of $Q^5$ 
by a proper side-pairing. Then the Euler characteristic of $M$ is zero.  
Although the $k$-faces of $Q^5$ are noncompact for all $k > 0$, 
we can compute the Euler characteristic of $M$ 
from the face-cycle decomposition of $M$, in the usual way, 
since the Euler characteristic of the link of each cusp of $M$ is zero. 
As a consistency check, we compute $\chi(M)$ 
in terms of the face-cycle decomposition of $M$,  
$$\chi(M) = 5-70+170-140+36-1 = 0.$$ 

By a computer search, we found that 
there are exactly 6,616,152 proper side-pairings of $Q^5$, 
all of whose elements are orientation preserving, 
of the form described in Lemma 1.   
These side-pairings fall into 55,168 equivalence classes 
under equivalence by a symmetry of $Q^5$. 
The classification of integral, congruence two, orientable,  
hyperbolic 5-manifolds of minimum volume is summarized 
in our next theorem. 

\begin{theorem}
There are, up to isometry, exactly 3607 
orientable hyperbolic space-forms $H^5/\Gamma$ 
where $\Gamma$ is a torsion-free subgroup of minimal index 
in the congruence two subgroup $\Gamma_2^5$ of the group 
$\Gamma^5$ of integral, positive, Lorentzian $6\times 6$ matrices. 
All these hyperbolic 5-manifolds have either 10 or 12 cusps, 
with only 26 of these manifolds having 12 cusps. 
All these hyperbolic 5-manifolds have volume $28\zeta(3)$. 
\end{theorem}
\begin{proof}
Each side-pairing of $Q^5$ induces an equivalence relation 
on the 90 ideal vertices of $Q^5$. The equivalence classes 
are cycles. The cycle of a large ideal vertex $v = \pm e_i+ e_6$ 
of $Q^5$ is either just itself or itself and its antipodal vertex, 
since an element of $K^5$ either fixes $v$ or maps $v$ to $-v$. 
It turns out that the 10 large ideal vertices of $Q^5$  
fall into either 5 or 6 cycles of the form $2,2,2,2,2$ or $1,1,2,2,2,2$. 
The remaining 80 ideal vertices of $Q^5$ fall into either 5 or 6 cycles 
of the form $16,16,16,16,16$ or $8,8,16,16,16,16$. 
The possible vertex cycle structures are $2,2,2,2,2,16,16,16,16,16$ 
and $1,1,2,2,2,2,8,8,16,16,16,16$. 
Thus all the integral, congruence two, orientable, hyperbolic 
5-manifolds of minimum volume have either 10 or 12 cusps. 
Only 26 of the manifolds have 12 cusps. 

Let $M$ be a 10-cusped manifold.  Its side-pairing has a vertex 
cycle structure $2,2,2,2,2,16,16,16,16,16$. 
Consider a maximum open cusp of $M$ of vertex cycle order 2 
represented by two antipodal large ideal vertices $\pm v$ of $Q^5$. 
This maximum open cusp is represented by the two open horoballs 
based at the ideal vertices $\pm v$ whose boundaries are tangent 
at the center $e_6$ of $Q^5$. 
Thus the boundary horosphere of this maximum open cusp is tangent to itself 
at the point represented by the origin $e_6$ of $H^5$. 
For example, if $v=e_1+e_6$, then the equations of the supporting 
hyperplanes of the horospheres are $\pm x_1-x_6=-1$ which 
have $e_6$ as their only common solution in $H^5$. 

The two horospheres based at $\pm v$ that are tangent at the origin 
also pass through 32 actual vertices of $Q^5$, and these 32 vertices 
represent another self-tangency point of the maximum cusp. 
For example, if $v=e_1+e_6$, then the 32  actual vertices are 
$(\pm 2,\pm 1,\pm 1,\pm 1,\pm 1,3)$. 
These two horospheres are also tangent to 48 centers of ideal square faces of $Q^5$, 
and these points represent 6 more self-tangency points of the maximum cusp. 
For example, if $v=e_1+e_6$, then the 48 centers of ideal square faces 
are the 48 points obtained from $(\pm 1,0,0,\pm 1,\pm 1,2)$ by 
permuting the middle four coordinates. 
Thus, the boundary of a maximum open cusp of $M$ of vertex cycle 
order 2 is tangent to itself at 8 points. 

The boundary of a maximum open cusp of $M$ of vertex cycle 
order 16 is also tangent to itself at 8 points. 
For example, if the cycle of ideal vertices consists of the 
16 ideal vertices $(0,\pm 1,\pm 1,\pm 1, \pm1, 2)$, 
the maximum open cusp is represented by 16 open horoballs 
based at the ideal vertices $(0,\pm 1,\pm 1,\pm 1, \pm1, 2)$ 
whose horosphere boundaries are supported by the hyperplanes whose 
equations are 
$$\pm x_2 \pm x_3 \pm x_4 \pm x_5 -2x_6 = -1.$$
Each of these horospheres is tangent to 4 adjacent horospheres.  
These $4\cdot 16/2= 32$ tangency points are the 32 
points obtained from $(0,0,\pm 1,\pm 1,\pm 1,2)$ by permuting 
the middle four coordinates. These 32 points represent four 
self-tangency points of the maximum cusp. 
Each horosphere passes through 8 actual vertices of $Q^5$. 
These $8\cdot 16 =128$ vertices are obtained from 
$(\pm 1,\pm 2,\pm 1,\pm 1,\pm 1,3)$ 
by permuting the middle four coordinates.  
These 128 vertices represent another 4 self-tangency 
points of the maximum cusp. 

All the self-tangency points of the maximum cusps of $M$ consist 
of only 16 points. Each of these 16 points is a self-tangency point 
of the boundary of exactly 5 of the maximum cusps. 
Thus $M$ has a set of 16 canonical points. 
The 16 canonical points of $M$ are represented by the 16 actual vertices 
of the polytope $P^5$. 

Let $M$ be a 12-cusped manifold.  Its side-pairing has a vertex 
cycle structure $1,1,2,2,2,2,8,8,16,16,16,16$. 
Consider a maximum open cusp of $M$ of vertex cycle order 1 
represented by a large ideal vertex $v$ of $Q^5$.    
This maximum open cusp is represented by the open horoball  
based at $v$ whose boundary passes through the centers 
of the 24 large sides of $Q^5$ that are not incident with $\pm v$. 
These 24 points represent 12 self-tangency points of the boundary 
of the maximum open cusp. 
For example, if $v = e_1 + e_6$, the maximum open cusp is represented 
by the open horoball based at $v$ whose horosphere boundary is supported 
by the hyperplane whose equation is $x_1-x_6 = -\sqrt{2}$. 
The 24 centers of large sides that are not incident with $\pm v$ 
are obtained from $(0,\pm\sqrt{2}/2,\pm\sqrt{2}/2,0,0,\sqrt{2})$ 
by permuting the middle four coordinates. 
These 24 points represent 12 self-tangency points of the maximum cusp. 
Thus, the boundary of a maximum open cusp of $M$ of vertex cycle order 1 
is tangent to itself in at least 12 points. 

The boundary of a maximum open cusp of $M$ of vertex cycle order 8 is also 
tangent to itself in at least 12 points. For example, if the cycle of ideal vertices 
consists of the 8 ideal vertices $(0,\pm 1,\pm 1,\pm 1,\pm 1, 2)$, with an even 
number of minus signs, the maximum open cusp is represented by 8 open horoballs, 
based at these ideal vertices, whose horosphere boundaries are supported by 
the hyperplanes whose equations are 
$$\pm x_2 \pm x_3 \pm x_4 \pm x_5 -2x_6 = -\sqrt{2},$$
with an even number of minus signs. 
Each of these horospheres are tangent to 6 adjacent horospheres. 
These $6\cdot 8/2 = 24$ tangency points are the 24 points obtained from 
$(0,\pm\sqrt{2}/2,\pm\sqrt{2}/2,0,0,\sqrt{2})$ by permuting the middle 
four coordinates. These 24 points are centers of large sides of $Q^5$ 
and represent 12 self-tangency points of the maximum cusp of $M$. 

The boundary of a maximum open cusp of $M$ of vertex cycle order 2 or 16 
is tangent to itself at 8 points as in the 10 cusp case. 
Thus, the cusps of vertex cycle order 1 or 8 are intrinsically different 
from the cusps of vertex cycle order 2 or 16. 
In fact, the volume of a maximum cusp of vertex cycle order 1 or 8 
is 16 whereas the volume of a maximum cusp of vertex cycle order 2 or 16 is 8. 
Thus, a 12-cusped manifold has 4 large maximum cusps and 8 small maximum cusps.  
Now an isometry between 12-cusped manifolds 
must map small maximum cusps to small maximum cusps, and so 
the self-tangency points of the boundaries of the small maximum cusps 
are canonical points of a 12-cusped manifold. 
There are a total of 16 canonical points of a 12-cusped manifold  
represented by the 16 actual vertices of the polytope $P^5$. 

Now let $M$ be either a 10- or 12-cusped manifold. 
Each actual vertex of $P^5$ is the vertex of a right-angled corner. 
This suggests a cut and paste operation on $Q^5$. 
Cut $Q^5$ along the 5 coordinate hyperplanes into 32 copies of 
the polytope $P^5$. By reassembling a copy of $Q^5$ around 
a different corner of $P^5$ than the origin 
according to the gluing pattern of the side-pairing of $Q^5$ that glues up to $M$, 
we get possibly 16 different side-pairings of $Q^5$ that glue up to $M$. 
We call such a cut and paste operation an {\it inside-out operation} on $M$. 

We have 55,168 equivalence classes of side-pairings of $Q^5$ 
under equivalence by a symmetry of $Q^5$.
Each of these side-pairings of $Q^5$ determines 16 side-pairings on $Q^5$ 
after inside-out operations on $Q^5$ that yield the same manifold. 
After comparing the new side-pairings with the old ones, up to symmetry of $Q^5$, 
the number of manifolds is reduced to 3607. 

Let $\phi:M \to M'$ be an isometry between two of these 10- or 12-cusped manifolds. 
Let $c$ and $c'$ be the point of $M$ and $M'$, respectively, 
represented by the center $e_6$. 
Now $c$ is a canonical point of $M$, and so $\phi(c)$ is a canonical point of $M'$. 
By applying an inside-out operation to $M'$ 
if $\phi(c)\neq c'$, we may assume that $\phi(c) = c'$. 
Let $\tilde\phi:H^5\to H^5$ be the lift of $\phi$ such that $\tilde\phi(e_6)=e_6$. 
Let $\Gamma$ and $\Gamma'$ be the discrete groups generated by the side-pairings 
of $Q^5$ that glue up $M$ and $M'$, respectively. 
Then $Q^5$ is the Dirichlet polytope for $\Gamma$ and $\Gamma'$ 
centered at $e_6$ by Lemma 1. 
As $\Gamma'=\tilde\phi\Gamma\tilde\phi\hbox{}^{-1}$, 
we have that $\tilde\phi$ maps the Dirichlet polytope for $\Gamma$ centered at $e_6$ 
to the Dirichlet polytope of $\Gamma'$ centered at $e_6$. 
Therefore $\tilde\phi$ is a symmetry of $Q^5$. 
As the 3607 manifolds have been classified up to symmetry of $Q^5$, 
we deduce that $\tilde\phi$ is a symmetry of the side-pairing of $Q^5$ for $M$, 
and so $M$ = $M'$. 
Thus the classification of these manifolds up to isometry is complete.

\end{proof} 

\section{Side-Pairing Coding} 

In this section, we describe the coding that we use to specify 
the side-pairings of $Q^5$ for the integral, congruence two, 
hyperbolic 5-manifolds of minimum volume. 
Reading this section is necessary only if the reader 
wants to reconstruct one of the manifolds in Tables 3,\ldots,6.

\begin{table}
$$\begin{array}{lllllll}
s_{1\p{0}}&=&(\p{-}1,\p{-}1,\p{-}0,\p{-}0,\p{-}0,\p{-}1) & \quad & s_{37}&=&(\p{-}0,\p{-}0,\p{-}0,\p{-}1,\p{-}1,\p{-}1)\\
s_{2\p{0}}&=&(-1,\p{-}1,\p{-}0,\p{-}0,\p{-}0,\p{-}1) & \quad & s_{38}&=&(\p{-}0,\p{-}0,\p{-}0,-1,\p{-}1,\p{-}1)\\
s_{3\p{0}}&=&(\p{-}1,-1,\p{-}0,\p{-}0,\p{-}0,\p{-}1) & \quad & s_{39}&=&(\p{-}0,\p{-}0,\p{-}0,\p{-}1,-1,\p{-}1)\\
s_{4\p{0}}&=&(-1,-1,\p{-}0,\p{-}0,\p{-}0,\p{-}1) & \quad & s_{40}&=&(\p{-}0,\p{-}0,\p{-}0,-1,-1,\p{-}1)\\
s_{5\p{0}}&=&(\p{-}1,\p{-}0,\p{-}1,\p{-}0,\p{-}0,\p{-}1) & \quad & s_{41}&=&(\p{-}1,\p{-}1,\p{-}1,\p{-}1,\p{-}1,\p{-}2)\\
s_{6\p{0}}&=&(-1,\p{-}0,\p{-}1,\p{-}0,\p{-}0,\p{-}1) & \quad & s_{42}&=&(-1,\p{-}1,\p{-}1,\p{-}1,\p{-}1,\p{-}2)\\
s_{7\p{0}}&=&(\p{-}1,\p{-}0,-1,\p{-}0,\p{-}0,\p{-}1) & \quad & s_{43}&=&(\p{-}1,-1,\p{-}1,\p{-}1,\p{-}1,\p{-}2)\\
s_{8\p{0}}&=&(-1,\p{-}0,-1,\p{-}0,\p{-}0,\p{-}1) & \quad & s_{44}&=&(-1,-1,\p{-}1,\p{-}1,\p{-}1,\p{-}2)\\
s_{9\p{0}}&=&(\p{-}0,\p{-}1,\p{-}1,\p{-}0,\p{-}0,\p{-}1) & \quad & s_{45}&=&(\p{-}1,\p{-}1,-1,\p{-}1,\p{-}1,\p{-}2)\\
s_{10}&=&(\p{-}0,-1,\p{-}1,\p{-}0,\p{-}0,\p{-}1) & \quad & s_{46}&=&(-1,\p{-}1,-1,\p{-}1,\p{-}1,\p{-}2)\\
s_{11}&=&(\p{-}0,\p{-}1,-1,\p{-}0,\p{-}0,\p{-}1) & \quad & s_{47}&=&(\p{-}1,-1,-1,\p{-}1,\p{-}1,\p{-}2)\\
s_{12}&=&(\p{-}0,-1,-1,\p{-}0,\p{-}0,\p{-}1) & \quad & s_{48}&=&(-1,-1,-1,\p{-}1,\p{-}1,\p{-}2)\\
s_{13}&=&(\p{-}1,\p{-}0,\p{-}0,\p{-}1,\p{-}0,\p{-}1) & \quad & s_{49}&=&(\p{-}1,\p{-}1,\p{-}1,-1,\p{-}1,\p{-}2)\\
s_{14}&=&(-1,\p{-}0,\p{-}0,\p{-}1,\p{-}0,\p{-}1) & \quad & s_{50}&=&(-1,\p{-}1,\p{-}1,-1,\p{-}1,\p{-}2)\\
s_{15}&=&(\p{-}1,\p{-}0,\p{-}0,-1,\p{-}0,\p{-}1) & \quad & s_{51}&=&(\p{-}1,-1,\p{-}1,-1,\p{-}1,\p{-}2)\\
s_{16}&=&(-1,\p{-}0,\p{-}0,-1,\p{-}0,\p{-}1) & \quad & s_{52}&=&(-1,-1,\p{-}1,-1,\p{-}1,\p{-}2)\\
s_{17}&=&(\p{-}0,\p{-}1,\p{-}0,\p{-}1,\p{-}0,\p{-}1) & \quad & s_{53}&=&(\p{-}1,\p{-}1,-1,-1,\p{-}1,\p{-}2)\\
s_{18}&=&(\p{-}0,-1,\p{-}0,\p{-}1,\p{-}0,\p{-}1) & \quad & s_{54}&=&(-1,\p{-}1,-1,-1,\p{-}1,\p{-}2)\\
s_{19}&=&(\p{-}0,\p{-}1,\p{-}0,-1,\p{-}0,\p{-}1) & \quad & s_{55}&=&(\p{-}1,-1,-1,-1,\p{-}1,\p{-}2)\\
s_{20}&=&(\p{-}0,-1,\p{-}0,-1,\p{-}0,\p{-}1) & \quad & s_{56}&=&(-1,-1,-1,-1,\p{-}1,\p{-}2)\\
s_{21}&=&(\p{-}0,\p{-}0,\p{-}1,\p{-}1,\p{-}0,\p{-}1) & \quad & s_{57}&=&(\p{-}1,\p{-}1,\p{-}1,\p{-}1,-1,\p{-}2)\\
s_{22}&=&(\p{-}0,\p{-}0,-1,\p{-}1,\p{-}0,\p{-}1) & \quad & s_{58}&=&(-1,\p{-}1,\p{-}1,\p{-}1,-1,\p{-}2)\\
s_{23}&=&(\p{-}0,\p{-}0,\p{-}1,-1,\p{-}0,\p{-}1) & \quad & s_{59}&=&(\p{-}1,-1,\p{-}1,\p{-}1,-1,\p{-}2)\\
s_{24}&=&(\p{-}0,\p{-}0,-1,-1,\p{-}0,\p{-}1) & \quad & s_{60}&=&(-1,-1,\p{-}1,\p{-}1,-1,\p{-}2)\\
s_{25}&=&(\p{-}1,\p{-}0,\p{-}0,\p{-}0,\p{-}1,\p{-}1) & \quad & s_{61}&=&(\p{-}1,\p{-}1,-1,\p{-}1,-1,\p{-}2)\\
s_{26}&=&(-1,\p{-}0,\p{-}0,\p{-}0,\p{-}1,\p{-}1) & \quad & s_{62}&=&(-1,\p{-}1,-1,\p{-}1,-1,\p{-}2)\\
s_{27}&=&(\p{-}1,\p{-}0,\p{-}0,\p{-}0,-1,\p{-}1) & \quad & s_{63}&=&(\p{-}1,-1,-1,\p{-}1,-1,\p{-}2)\\
s_{28}&=&(-1,\p{-}0,\p{-}0,\p{-}0,-1,\p{-}1) & \quad & s_{64}&=&(-1,-1,-1,\p{-}1,-1,\p{-}2)\\
s_{29}&=&(\p{-}0,\p{-}1,\p{-}0,\p{-}0,\p{-}1,\p{-}1) & \quad & s_{65}&=&(\p{-}1,\p{-}1,\p{-}1,-1,-1,\p{-}2)\\
s_{30}&=&(\p{-}0,-1,\p{-}0,\p{-}0,\p{-}1,\p{-}1) & \quad & s_{66}&=&(-1,\p{-}1,\p{-}1,-1,-1,\p{-}2)\\
s_{31}&=&(\p{-}0,\p{-}1,\p{-}0,\p{-}0,-1,\p{-}1) & \quad & s_{67}&=&(\p{-}1,-1,\p{-}1,-1,-1,\p{-}2)\\
s_{32}&=&(\p{-}0,-1,\p{-}0,\p{-}0,-1,\p{-}1) & \quad & s_{68}&=&(-1,-1,\p{-}1,-1,-1,\p{-}2)\\
s_{33}&=&(\p{-}0,\p{-}0,\p{-}1,\p{-}0,\p{-}1,\p{-}1) & \quad & s_{69}&=&(\p{-}1,\p{-}1,-1,-1,-1,\p{-}2)\\
s_{34}&=&(\p{-}0,\p{-}0,-1,\p{-}0,\p{-}1,\p{-}1) & \quad & s_{70}&=&(-1,\p{-}1,-1,-1,-1,\p{-}2)\\
s_{35}&=&(\p{-}0,\p{-}0,\p{-}1,\p{-}0,-1,\p{-}1) & \quad & s_{71}&=&(\p{-}1,-1,-1,-1,-1,\p{-}2)\\
s_{36}&=&(\p{-}0,\p{-}0,-1,\p{-}0,-1,\p{-}1) & \quad & s_{72}&=&(-1,-1,-1,-1,-1,\p{-}2)
\end{array}$$
\caption{The Lorentz normal vectors of the 72 sides of $Q^5$}
\end{table}

We encode an element $\diag(a_{1},a_{2},a_{3},a_{4},a_{5},1)$ of $K^{5}$, where 
each $a_{i}=\pm 1$, by
$$\sum_{i=1}^{5}\frac{(1-a_{i})}{2}2^{i-1}$$
expressed as a single base 32 digit $0, \ldots, 9,A=10, \ldots, V=31$.

Table 1 gives the Lorentz normal vector $s_i$ of the $i$th side $S_i$ of $Q^5$ 
for $i = 1,\ldots 72$.  
The side-pairing map that maps the side $S'_i$ to the side $S_{i}$ 
is of the form $r_{i}k_{i}$ where 
$r_{i}$ is the reflection in side $S_{i}$ and $k_{i}$ is an element of $K^{5}$  
such that $S'_i = k_iS_i$.   
We have necessarily in these manifolds that the large sides pair in 
groups of four so that 
$k_{4i-3}=k_{4i-2}=k_{4i-1}=k_{4i}$
for $i=1,2,\ldots 10$ and 
then the remaining 32 small sides all pair in a similar fashion 
with $k_{41}=k_{42}=\cdots=k_{72}$.
Thus it suffices to take 
$$k_{1}k_{5}k_{9}k_{13}k_{17}k_{21}k_{25}k_{29}k_{33}k_{37}k_{41}$$ 
as an eleven digit base 32 code for the side-pairing.

By permuting coordinates, a given manifold may have several different 
encodings.  Moreover there is an inside-out operation that takes the 
$Q^{5}$ fundamental domain apart and rearranges the pieces to give 
another side-pairing of $Q^{5}$ for the same manifold.  We take only 
one representative coding for each manifold and count the number of 
symmetries of the manifold that give the same side-pairing (these 
are all of the possible symmetries of the manifold).  Thus the 
last digit of the code represents how the small sides pair, and by a 
symmetry we may as well take $k_{41}$ coded by $1$, $7$, or $V=31$, 
i.e., the small sides map with $k_{i}=\diag(-1,1,1,1,1,1)$, 
$\diag(-1,-1,-1,1,1,1)$, or $\diag(-1,-1,-1,-1,-1,1)$. 

Finally, it is worth noting that the form of the side-pairing maps described above implies 
that all the torsion-free subgroups of $\Gamma^5_2$ of minimal index are normal.

\section{Tables of Manifolds} 
Tables 3-6 list side-pairings and isometric invariants of 
some examples of integral, congruence two, hyperbolic 5-manifolds of minimum volume. 
In each table, $N$ is the row number in the complete list 
of all integral, congruence two, hyperbolic 5-manifolds of minimum volume. 
The 12-cusped manifolds have row numbers $1,\ldots,26$. 
The column headed by $SP$ lists the side-pairing for the manifold 
in a coded form that is explained in the previous section. 
The column headed by $S$ lists the number of symmetries of the manifolds. 
All the manifolds have a subgroup of symmetries corresponding to $K^5$. 
Therefore, the number of symmetries is a multiple of 32. 
The possible orders are $32, 64, 128, 160,192,256,384,512,1536$. 
We ordered the 10- and 12-cusped manifolds so that the number of 
symmetries $S$ is nonincreasing. 

The column headed by $H_i$ lists the $i$-th homology groups of the manifolds 
with the three numbers $a, b, c$ representing $\integers^a\oplus\integers^b_2\oplus\integers^c_4$ 
and the single entry $a$ representing $\integers^a$. 
All the 12-cusped manifolds have their fourth homology group isomorphic to $\integers^{11}$. 
All the 10-cusped manifolds have their fourth homology group isomorphic to $\integers^9$. 
The column headed by $LT$ lists the link types of the links of the cusps of the manifolds. 
In Table 3, the four large cusp link types are listed last. 

There are 27 orientable, closed, flat 4-manifolds up to homeomorphism
according to J. Hillman [2] and R. Levine [4].  
We denote these 27 manifolds by $O^4_1,\ldots,O^4_{27}$ in the order 
of Table D in Levine's thesis [4].  We also denote manifolds 
$O^4_1,\ldots,O^4_{26}$ by the letters $A,\ldots, Z$. 
Generators for $\pi_1(O^4_i)$, as a group of affine transformations of $\realnos^4$, 
are given by Levine in Tables C and D of his thesis [4]. 
The first homology and holonomy groups 
of the closed orientable flat 4-manifolds were computed by Levine 
and are listed in Table D of his thesis [4]. 

The fundamental group of the link of a cusp of any one of our 5-manifolds is  
represented by a parabolic subgroup of $\Gamma^5_2$, and so is conjugate to a subgroup 
of a stabilizer of an ideal vertex of $P^5$.  
The link of an ideal vertex of $P^5$ is a Euclidean 4-cube, 
and so the stabilizer of an ideal vertex of $P^5$ is the reflection 
group in the sides of a Euclidean 4-cube. 
The holonomy group of the reflection group in the sides 
of a Euclidean 4-cube with one vertex at the origin and sides on 
the four coordinate hyperplanes of $\realnos^4$ is an elementary 2-group of rank 4 
generated by the reflections in the four coordinate hyperplanes of $\realnos^4$. 
Thus, the holonomy groups of the links of the cusps of our 5-manifolds 
are elementary 2-groups. 
These groups have rank at most two by Levine's computation of holonomy groups. 
Thus, the holonomy groups of the links of the cusps of our 5-manifolds 
are elementary 2-groups of rank 0, 1, or 2. 
The closed orientable flat 4-manifolds, with an elementary 2 holonomy group, 
are the manifolds $A, \ldots, L$.  

We denote the six orientable, closed, flat 3-manifolds by $O^3_1,\ldots, O^3_6$ 
in the order given by J. Wolf [13].  
Each of the manifolds $O^4_i$ for $i=1,\ldots, 27$ is a mapping torus over $O^3_j$ 
for some $j = 1,\ldots 6$; moreover, if the first Betti number of $O^4_i$ is 1, 
then $O^4_i$ is a mapping torus over $O^3_j$ for just one $j$. 
See Hillman [2] for a discussion. 
We denote a nondirect product mapping torus over $O^3_j$ by $O^3_j\rtimes S^1$. 
A mapping torus decomposition of $O^4_i$, either $O^3_j\times S^1$ or $O^3_j\rtimes S^1$, 
corresponds, respectively, to either a direct product $\pi(O^3_j)\times\integers$ 
or a semidirect product $\pi(O^3_j)\rtimes\integers$ decomposition of $\pi_1(O^4_i)$. 

Table 2 lists the homeomorphism type, first homology group, mapping torus type, 
and holonomy group of the manifolds $A,\ldots, L$.  
The mapping torus types of the manifolds $A,\ldots, L$ were determined by expressing 
$\pi_1(O^4_i)$ as either a direct product $\pi(O^3_j)\times\integers$ 
or a semidirect product $\pi(O^3_j)\rtimes\integers$. 

The manifolds $A$ and $B$ are determined by their first homology groups. 
The manifold $C$ is determined by its first homology and holonomy groups. 
The manifolds $D$ and $E$ are determined by their first homology groups.  
The manifolds $F$ and $G$ have the same first homology and holonomy groups and mapping torus type. 
We distinguished these manifolds by their fundamental groups. 
The manifold $G$ does not occur as a link type of one of our 5-manifolds.  
The manifold $H$ is determined by its first homology and holonomy groups and mapping torus type. 
The manifolds $I$ and $J$ have the same first homology and holonomy groups and mapping torus type. 
We distinguished these manifolds by their fundamental groups. 
The manifold $K$ is determined by its first homology and holonomy groups and mapping torus type.
The manifold $L$ is determined by its first homology and holonomy groups. 
The manifold $L$ does not occur as a link type of one of our 5-manifolds.  

The coordinate hyperplane cross-sections of the links of the cusps of our 5-manifolds 
are totally geodesic hypersurfaces.  The mapping torus types of the links of the cusps 
of our 5-manifolds were determined by finding a connected orientable cross-section 
with a connected complement. Such a cross-section determines a mapping torus decomposition 
of the link of a cusp, since the link of each large ideal vertex of $Q^5$ is a 4-cube 
and by an inside-out operation on $Q^5$ we can represent any cusp by a cycle of large 
ideal vertices of $Q^5$.

\begin{table} 
\begin{center}
\begin{tabular}{lllll}
\rule{0pt}{15pt}$LT$&$HT$&$H_{1}$&$MT$&$HG$\\

\rule{0pt}{15pt}A&$O^4_1$&$\mathbb{Z}^4$                    &$O^3_1\times S^1$&1\\
\rule{0pt}{15pt}B&$O^4_2$&$\mathbb{Z}^2\oplus\mathbb{Z}_2^2$&$O^3_2\times S^1$&$\integers_2$\\
\rule{0pt}{15pt}C&$O^4_3$&$\mathbb{Z}^2\oplus\mathbb{Z}_2$  &$O^3_2\rtimes S^1$&$\integers_2$\\
\rule{0pt}{15pt}D&$O^4_4$&$\mathbb{Z}\oplus\mathbb{Z}_4^2$  &$O^3_6\times S^1$&$\integers^2_2$\\
\rule{0pt}{15pt}E&$O^4_5$&$\mathbb{Z}\oplus\mathbb{Z}_2^3$  &$O^3_2\rtimes S^1$&$\integers^2_2$\\
\rule{0pt}{15pt}F&$O^4_6$&$\mathbb{Z}\oplus\mathbb{Z}_2\oplus\mathbb{Z}_4$&$O^3_2\rtimes S^1$&$\integers^2_2$\\
\rule{0pt}{15pt}G&$O^4_7$&$\mathbb{Z}\oplus\mathbb{Z}_2\oplus\mathbb{Z}_4$&$O^3_2\rtimes S^1$&$\integers^2_2$\\
\rule{0pt}{15pt}H&$O^4_8$&$\mathbb{Z}\oplus\mathbb{Z}_2\oplus\mathbb{Z}_4$&$O^3_6\rtimes S^1$&$\integers^2_2$\\
\rule{0pt}{15pt}I&$O^4_9$&$\mathbb{Z}\oplus\mathbb{Z}_2^2$  &$O^3_6\rtimes S^1$&$\integers^2_2$\\
\rule{0pt}{15pt}J&$O^4_{10}$&$\mathbb{Z}\oplus\mathbb{Z}_2^2$&$O^3_6\rtimes S^1$&$\integers^2_2$\\
\rule{0pt}{15pt}K&$O^4_{11}$&$\mathbb{Z}\oplus\mathbb{Z}_2^2$&$O^3_2\rtimes S^1$&$\integers^2_2$\\
\rule{0pt}{15pt}L&$O^4_{12}$&$\mathbb{Z}\oplus\mathbb{Z}_4$    &$O^3_2\rtimes S^1$&$\integers^2_2$\\
\end{tabular}
\end{center}

\vspace{.1in}
\caption{The closed orientable flat 4-manifolds with elementary 2 holonomy}
\end{table}

\begin{table}
\begin{displaymath}
\begin{array}{rcrrrrrrrrl}
N&SP&S&\ H_{1}&&&\ H_{2}&&&\ H_{3}&LT\\
1&{\tt 24B8DPGPDB1}&1536&4&7&0&12&4&6&20&{\tt JJJJJJJJAAAA}\\
2&{\tt PVSE8BBGEQ7}&512&0&10&1&12&10&0&24&{\tt BBBBCCCCBBBB}\\
3&{\tt 1B4EDQVLEB7}&256&2&9&0&4&8&8&14&{\tt HHHHJJJJBBBB}\\
4&{\tt 1B4EDQQGEB7}&128&2&9&0&4&8&8&14&{\tt HHJJJJJJBBBB}\\
5&{\tt 1BLEPVQGEB7}&128&1&10&0&6&8&5&17&{\tt CCJJJJKKBBBB}\\
6&{\tt 1B4ESQQLEV7}&128&1&10&0&2&12&6&13&{\tt BBJJJJKKEEEE}\\
7&{\tt 1B4EDMJGEB7}&64&1&10&0&6&6&6&17&{\tt HHHHKKKKBBBB}\\
8&{\tt 1B4EDMQPEB7}&64&1&10&0&4&8&6&15&{\tt HHHJKKKKEEEE}\\
9&{\tt 1EDQ8BBGME7}&64&1&10&0&2&9&8&13&{\tt JJJJKKKKBBEE}\\
10&{\tt 1B4EDQJPEB7}&64&1&10&0&0&10&9&11&{\tt HHHHJJJJEEEE}\\
11&{\tt 1EDB8MMLJE7}&64&0&11&0&4&10&4&16&{\tt BCHHKKKKEEEE}\\
12&{\tt 1B4EPMJGEB7}&64&0&11&0&4&8&6&16&{\tt HHHHKKKKBBBB}\\
13&{\tt 1BSEPVJGEB7}&64&0&11&0&4&8&5&16&{\tt CCHHHHKKEEEE}\\
14&{\tt 1B4EDMMLEJ7}&64&0&11&0&4&7&6&16&{\tt HHHHKKKKBBEE}\\
15&{\tt 1BLEPMJGEB7}&64&0&11&0&4&6&7&16&{\tt HHHHKKKKBBBB}\\
16&{\tt 1BSEPMQGEB7}&64&0&11&0&3&14&3&15&{\tt BBCJJKKKEEEE}\\
17&{\tt 1B4ESQMLEB7}&64&0&11&0&1&8&8&13&{\tt HHHHHJKKEEEE}\\
18&{\tt 1EDJ8MBLJE7}&64&0&11&0&1&8&8&13&{\tt HJKKKKKKEEEE}\\
19&{\tt 1B4EDMJPEB7}&32&1&10&0&3&7&8&14&{\tt HHHHJJKKBBEE}\\
20&{\tt 1B4EDMMSEB7}&32&1&10&0&3&7&8&14&{\tt HHHJJJKKBBEE}\\
21&{\tt 1BSEPMJGEB7}&32&0&11&0&7&10&2&19&{\tt BBCCHHKKBBEE}\\
22&{\tt 174EDQQPEJV}&32&0&11&0&3&10&5&15&{\tt CHHHHKKKEEEE}\\
23&{\tt 1BLEPMQGEB7}&32&0&11&0&3&9&6&15&{\tt CHHJJKKKBBEE}\\
24&{\tt 1B4EDQQPEV7}&32&0&11&0&2&8&7&14&{\tt CHHHJJJKEEEE}\\
25&{\tt 174EDQQGEJV}&32&0&11&0&2&7&8&14&{\tt HHJJKKKKBBEE}\\
26&{\tt 1B4EPMMLQBV}&32&0&11&0&1&8&8&13&{\tt HHJJKKKKEEEE} 
\end{array}
\end{displaymath}
\caption{The orientable, 12-cusped, minimum volume, integral, 
congruence 2, hyperbolic 5-manifolds}
\end{table}

\begin{table}
\begin{displaymath}
\begin{array}{rcrrrrrrrrl}
N&SP&S&H_{1}&&&H_{2}&&&H_{3}&LT\\
27&{\tt 2B7JB47JG81}&1536&4&7&0&6&7&9&12&{\tt AADDDDDDDD}\\
28&{\tt 24B8DPPGDB1}&384&4&7&0&6&6&8&12&{\tt AADDJJJJJJ}\\
29&{\tt 27BBJ4J7G81}&384&4&7&0&6&5&9&12&{\tt AADDDDDDJJ}\\
30&{\tt 24BBJ7J7G81}&256&4&7&0&6&6&8&12&{\tt AADDDDJJJJ}\\
31&{\tt 1DJEDBS2EV7}&256&2&9&0&10&8&1&18&{\tt AACCCCJJJJ}\\
32&{\tt 24BJS8BJSG1}&256&2&9&0&10&7&2&18&{\tt AACCCCHHHH}\\
33&{\tt 27BS74BSG81}&256&2&9&0&6&9&4&14&{\tt CCCCCCDDDD}\\
34&{\tt 1EDDJBV2SE7}&256&2&9&0&6&8&4&14&{\tt CCCCCCIIII}\\
35&{\tt DES12BJEDV7}&256&2&9&0&6&4&8&14&{\tt AAIIIIJJJJ}\\
36&{\tt 1ESDEBV2DJ7}&256&2&9&0&2&9&8&10&{\tt CCDDDDIIII}\\
37&{\tt 1DEESBJ2VD7}&256&2&9&0&2&6&10&10&{\tt CCDDDDHHHH}\\
38&{\tt BDV8DBGQSP7}&256&0&11&0&10&6&2&20&{\tt AABBBBKKKK}\\
39&{\tt BDV824GLJP7}&256&0&10&1&10&5&3&20&{\tt AABBBBHHHH}\\
40&{\tt DJQBJMS2487}&192&0&11&0&3&4&8&13&{\tt ADIIIIIIJJ}\\
41&{\tt DBQQMLMLDBV}&160&0&11&0&5&6&4&15&{\tt FFFFFKKKKK}\\
42&{\tt 14EPB8L2MGV}&160&0&11&0&0&6&9&10&{\tt EEEEEHHHHH}\\
43&{\tt DMQSJ7BLPEV}&160&0&11&0&0&6&9&10&{\tt IIIIIIIIII}\\
44&{\tt 1DJEDBJ2EG7}&128&2&9&0&10&6&2&18&{\tt AACCCCHHJJ}\\
45&{\tt 1DJEDBG2EJ7}&128&2&9&0&6&10&3&14&{\tt CCCCCCDDJJ}\\
46&{\tt 27BV74BVG81}&128&2&9&0&6&8&5&14&{\tt CCDDDDFFFF}\\
47&{\tt 27BS74BGS81}&128&2&9&0&6&8&4&14&{\tt BBCCCCDDHH}\\
48&{\tt 1EDDJBG2JE7}&128&2&9&0&6&8&4&14&{\tt BBCCCCDDII}\\
49&{\tt 1EDDJBJ2GE7}&128&2&9&0&6&8&4&14&{\tt BBCCCCHHII}\\
50&{\tt 1DJEDBV2ES7}&128&2&9&0&6&8&4&14&{\tt CCCCCCIIJJ}\\
51&{\tt 1DVEDBG2EV7}&128&2&9&0&6&7&6&14&{\tt AADDFFFFII}\\
52&{\tt 27BBS4S7G81}&128&2&9&0&6&7&5&14&{\tt BBDDFFFFII}\\
53&{\tt 1EDDVBV2GE7}&128&2&9&0&6&7&5&14&{\tt BBFFFFIIJJ}\\
54&{\tt 24BBJSJSG81}&128&2&9&0&6&6&6&14&{\tt AAFFFFHHII}\\
55&{\tt 1DVEDBV2EG7}&128&2&9&0&6&6&6&14&{\tt AAFFFFIIII}\\
56&{\tt 1EDDVBG2VE7}&128&2&9&0&6&6&6&14&{\tt CCDDFFFFJJ} 
\end{array}
\end{displaymath}
\caption{The first 60 orientable, 10-cusped, minimum volume, 
integral, congruence 2, hyperbolic 5-manifolds}
\end{table}

\begin{table}
\begin{displaymath}
\begin{array}{rcrrrrrrrrl}
N&SP&S&H_{1}&&&H_{2}&&&H_{3}&LT\\
57&{\tt 2B4JBDDJG81}&128&2&9&0&6&5&7&14&{\tt AADDFFFFHH}\\
58&{\tt DES12BSEDG7}&128&2&9&0&6&4&8&14&{\tt AAIIJJJJJJ}\\
59&{\tt 1EJDEBG2DJ7}&128&2&9&0&2&10&7&10&{\tt CCDDFFFFJJ}\\
60&{\tt 1DVEDBS2EJ7}&128&2&9&0&2&10&7&10&{\tt CCFFFFIIJJ}\\
61&{\tt 1ESDEBG2DS7}&128&2&9&0&2&8&9&10&{\tt CCDDDDDDII}\\
62&{\tt 1DEEJBG2JD7}&128&2&9&0&2&8&8&10&{\tt BBDDFFFFHH}\\
63&{\tt 1DEEJBJ2GD7}&128&2&9&0&2&8&8&10&{\tt BBFFFFHHHH}\\
64&{\tt 1EJDEBJ2DG7}&128&2&9&0&2&8&8&10&{\tt BBFFFFHHJJ}\\
65&{\tt 1EDDVBS2JE7}&128&2&9&0&2&8&8&10&{\tt BBFFFFJJJJ}\\
66&{\tt 1DEEJBV2SD7}&128&2&9&0&2&8&8&10&{\tt CCFFFFHHII}\\
67&{\tt 1DEEJBS2VD7}&128&2&9&0&2&8&8&10&{\tt CCFFFFHHJJ}\\
68&{\tt 1DEESBV2JD7}&128&2&9&0&2&7&9&10&{\tt BBDDDDHHII}\\
69&{\tt 1DEESBS2GD7}&128&2&9&0&2&7&9&10&{\tt BBDDDDHHJJ}\\
70&{\tt 1ESDEBS2DG7}&128&2&9&0&2&7&9&10&{\tt BBDDDDIIJJ}\\
71&{\tt 1B4EDQGQEB7}&128&2&9&0&2&6&10&10&{\tt BBDDHHJJJJ}\\
72&{\tt DES12BGEDS7}&128&2&9&0&2&6&10&10&{\tt BBDDIIJJJJ}\\
73&{\tt 1B4EDQLVEB7}&128&2&9&0&2&6&10&10&{\tt BBHHIIJJJJ}\\
74&{\tt 1DEESBG2SD7}&128&2&9&0&2&6&10&10&{\tt CCDDDDDDHH}\\
75&{\tt 1SV8MBBLEV7}&128&1&10&0&8&7&2&17&{\tt BBCCCCCCKK}\\
76&{\tt BDQED8PE4S7}&128&1&10&0&8&7&2&17&{\tt BBCCCCCCKK}\\
77&{\tt BDVE8BG2EV7}&128&1&10&0&4&8&5&13&{\tt BBCCFFFFII}\\
78&{\tt BDVE8BV2EG7}&128&1&10&0&4&8&5&13&{\tt BBCCFFFFII}\\
79&{\tt 14BBMEJGM8V}&128&0&11&0&10&8&1&20&{\tt AABBBCHHHH}\\
80&{\tt 1477QEJGQ8V}&128&0&11&0&10&8&0&20&{\tt BBBBBCKKKK}\\
81&{\tt BLQQ2S1L487}&128&0&11&0&10&6&2&20&{\tt AABBCCHHHH}\\
82&{\tt 1V48BMBMGV7}&128&0&11&0&8&6&2&18&{\tt BCCCCCHHHH}\\
83&{\tt BDVESQGLJV7}&128&0&11&0&6&12&1&16&{\tt BBBBCCJJJJ}\\
84&{\tt 1B48MSPBMG7}&128&0&11&0&6&12&1&16&{\tt BCEEEEKKKK}\\
85&{\tt 1V4JDMSVPG7}&128&0&11&0&6&11&2&16&{\tt BBBCCCKKKK}\\
86&{\tt BDQESPMVLJ7}&128&0&11&0&2&8&6&12&{\tt BBIIIIJJJJ} 
\end{array}
\end{displaymath}
\caption{The first 60 orientable, 10-cusped, minimum volume, 
integral, congruence 2, hyperbolic 5-manifolds (cont.)}
\end{table}

\begin{table}
\begin{displaymath}
\begin{array}{rcrrrrrrrrl}
N&SP&S&H_{1}&&&H_{2}&&&H_{3}&LT\\
100&{\tt 1BS8JMG24B7}&64&2&9&0&4&6&7&12&{\tt BBBCFFHHII}\\
200&{\tt 14EEPBSBLD7}&64&1&10&0&2&6&8&11&{\tt CCEEHHHHJJ}\\
300&{\tt 1B48BDQGVP7}&32&2&9&0&8&7&4&16&{\tt AABCEHIKKK}\\
400&{\tt 14B8JLJGDB7}&32&2&9&0&5&7&6&13&{\tt ABCEEHIKKK}\\
500&{\tt 14BB2LPQD87}&32&2&9&0&3&6&8&11&{\tt BBCDFHHIJK}\\
600&{\tt 1BSSE4V2EG7}&32&1&10&0&6&6&4&15&{\tt BCCDFHKKKK}\\
700&{\tt 14QDVSMED87}&32&1&10&0&5&7&5&14&{\tt BBCCEEHIKK}\\
800&{\tt 1BEEP4S2JB7}&32&1&10&0&5&6&5&14&{\tt CCCFFIJKKK}\\
900&{\tt 1BL82VLVQB7}&32&1&10&0&4&8&5&13&{\tt BCCCDFFIIK}\\
1000&{\tt 1427D8GQDSV}&32&1&10&0&4&6&6&13&{\tt BBCCFHHIIK}\\
1100&{\tt 17DDMEMEQ8V}&32&1&10&0&4&6&6&13&{\tt BCCDEEIIKK}\\
1200&{\tt 1DE82PQPSB7}&32&1&10&0&4&6&6&13&{\tt CCCCEEHIIJ}\\
1300&{\tt 14EESB7QMGV}&32&1&10&0&4&5&7&13&{\tt BCCEEHIIKK}\\
1400&{\tt 1DJ8V4SMLP7}&32&1&10&0&3&7&7&12&{\tt BCCEIIKKKK}\\
1500&{\tt 14Q8JDQGVB7}&32&1&10&0&3&6&7&12&{\tt BBFHHHHIIK}\\
1600&{\tt 1BM8BEPJSQ7}&32&1&10&0&3&6&7&12&{\tt BCDEHHIIJK}\\
1700&{\tt 14B7E8M7LQV}&32&1&10&0&3&6&7&12&{\tt CCCEEIIIJK}\\
1800&{\tt 14EJDBBLVE7}&32&1&10&0&2&8&7&11&{\tt BBDEHHIIIJ}\\
1900&{\tt 17DE74QEJPV}&32&1&10&0&2&6&8&11&{\tt BCDEFFIIKK}\\
2000&{\tt 17BJSES7L8V}&32&1&10&0&2&6&8&11&{\tt BDEFHHHIIK}\\
2100&{\tt 14EESBLBGD7}&32&1&10&0&2&6&8&11&{\tt CDEFFHHIKK}\\
2200&{\tt 1DQ8JMS24B7}&32&1&10&0&1&6&9&10&{\tt BEEHIIIIKK}\\
2300&{\tt 14QDV8VLGS7}&32&0&11&0&5&7&4&15&{\tt BCEEFHHIKK}\\
2400&{\tt 1BSEJMJ2ES7}&32&0&11&0&4&7&5&14&{\tt BCCFFFFHHJ}\\
2500&{\tt 14BPM8BPMQV}&32&0&11&0&4&6&5&14&{\tt CCCDFFHIIK}\\
2600&{\tt 14QDE8MPVQ7}&32&0&11&0&3&8&5&13&{\tt BCCFFHIJKK}\\
2700&{\tt 14EEJLLBGD7}&32&0&11&0&3&6&6&13&{\tt BCCFFHHIKK}\\
2800&{\tt 14BSJDJEMG7}&32&0&11&0&3&6&6&13&{\tt CCCEFHIIIK}\\
2900&{\tt 17BQD4MSQBV}&32&0&11&0&3&6&6&13&{\tt CEEFHHIIKK}\\
3000&{\tt 14EBMS7SPGV}&32&0&11&0&2&6&7&12&{\tt BBDFHIIIJK}\\
3100&{\tt 1DV82BQLEP7}&32&0&11&0&2&6&7&12&{\tt BCFHIIIKKK}\\
3200&{\tt 1BVEDMV2PG7}&32&0&11&0&2&6&7&12&{\tt CCFFFHIJJK}\\
3300&{\tt 1B7JQ4SJPEV}&32&0&11&0&2&6&7&12&{\tt CFFHHIIIII}\\
3400&{\tt 14EJPBP2GM7}&32&0&11&0&1&6&8&11&{\tt CDDDFFHHIJ}\\
3500&{\tt 1BMDBLMPQ8V}&32&0&11&0&1&6&8&11&{\tt CFFFHHIIIK}\\
3600&{\tt 17BQM4BQLPV}&32&0&11&0&0&6&9&10&{\tt FFFFHHHHIK} 
\end{array}
\end{displaymath}
\caption{A selection of the other orientable, 10-cusped, minimum volume, 
integral, congruence 2, hyperbolic 5-manifolds}
\end{table}

\section{A Non-Orientable Example}

Our first example of a minimal volume, integral, congruence two, hyperbolic 
5-manifold was a nonorientable manifold $M$ constructed by Steven Tschantz 
in 1993 by gluing together the sides of the polytope $Q^5$ by hand. 
The side-pairing code for $M$ is {\tt 2549A81IKGV}. 
The manifold $M$ is the nonorientable hyperbolic 5-manifold of volume $28\zeta(3)$ 
mentioned in Ratcliffe [8].  
The manifold $M$ has $H_1(M) = \integers^5\oplus\integers^6_2$, 
$H_2(M) = \integers^{10}\oplus\integers_2\oplus\integers^9_4$, 
$H_3(M) = \integers^6\oplus\integers^9_2$, and $H_4(M)=0$.  
The manifold $M$ has 10 cusps all of which are nonorientable. 
Five of the links of the cusps have their first homology group equal to 
$\mathbb{Z}^{3}\oplus\mathbb{Z}_{2}$ and five have their 
first homology group equal to $\mathbb{Z}^{2}\oplus\mathbb{Z}_{2}$. 
The manifold $M$ has at least 160 symmetries.  

In our paper [10], we showed that   
there are 3 orientable and 10 nonorientable minimal volume, 
integral, congruence two, hyperbolic 3-manifolds and  
there are 22 orientable and 1149 nonorientable minimal volume, 
integral, congruence two, hyperbolic 4-manifolds. 
In this paper, we showed that there are 3607 orientable 
minimal volume, integral, congruence two, hyperbolic 5-manifolds, and 
so we expect that there are hundreds of thousands of nonorientable,  
minimal volume, integral, congruence two, hyperbolic 5-manifolds. 
These nonorientable manifolds can be classified in the same way 
as the orientable manifolds, 
but the classification is beyond the reach of present computer technology.

\section{A Small Hyperbolic 5-Manifold}

In this section, we describe a hyperbolic 5-manifold, 
of volume $7\zeta(3)/4$, that is obtained as a quotient space of one 
of our hyperbolic 5-manifolds.  
Let $M$ be manifold number 27 which is the 10-cusped hyperbolic 5-manifold 
with the largest symmetry group listed first in Table 4. 
The group of symmetries of $M$ has a subgroup 
of 16 orientation preserving symmetries that acts freely on $M$. 
This group of order 16 is generated by two elements $\alpha$ and $\beta$ 
subject to the relations $\alpha^8 =1$, $\beta^2=1$, $\beta\alpha\beta=\alpha^3$. 
Lifts to $H^5$ of the isometries $\alpha$ and $\beta$ of $M$  
are represented by the Lorentzian matrices 

$$\left(\begin{array}{cccccc}
     \p{-}1 & 0 &\phantom{-}0 & \p{-}1 & 0 &     -1 \\
     \p{-}0 & 0 &\phantom{-}0 & \p{-}0 & 1 & \p{-}0 \\
         -1 & 0 &          -1 & \p{-}0 & 0 & \p{-}1 \\
     \p{-}0 & 1 &\phantom{-}0 & \p{-}0 & 0 & \p{-}0 \\
     \p{-}0 & 0 &          -1 &     -1 & 0 & \p{-}1 \\
		 -1 & 0 &          -1 &     -1 & 0 & \p{-}2 
        \end{array} \right),
\left(\begin{array}{cccccc}
          0 & 1 & \p{-}0 & \p{-}0 & \p{-}0 & \p{-}0 \\
          1 & 0 & \p{-}0 & \p{-}0 & \p{-}0 & \p{-}0 \\
          0 & 0 & \p{-}1 & \p{-}1 & \p{-}0 &     -1 \\
		  0 & 0 &     -1 & \p{-}0 &     -1 & \p{-}1 \\
		  0 & 0 & \p{-}0 &     -1 &     -1 & \p{-}1 \\
		  0 & 0 &     -1 &     -1 &     -1 & \p{-}2 
        \end{array} \right).$$

\vspace{.1in} 
Let $N$ be the quotient manifold under the action of this group of order 16. 
Then $N$ is an orientable hyperbolic 5-manifold with two cusps. 
The cusps of $N$ have type $D$ and $P$.  
We have $H_1(P)=\integers\oplus\integers_4$. 
The mapping torus type of $P$ is $O^3_6\rtimes S^1$. 
The holonomy group of $P$ is a dihedral group of order eight. 
The manifold $P$ is determined by its first homology group and its mapping torus type. 
The cusp of $N$ of type $D$ is covered by the two cusps of $M$ of type $A$ 
while the cusp of $N$ of type $P$ is covered by the 8 cusps of $M$ of type $D$.  
The manifold $N$ has $H_1(N)=\integers\oplus\integers_4$, 
$H_2(N)=\integers^2_4$, $H_3(N)=\integers$, $H_4(N)=\integers$. 

The polytope $Q^5$ is subdivided into 32 copies of the polytope $P^5$ 
by the five coordinate hyperplanes $x_i=0$ for $i=1,\ldots 5$, 
and so the manifold $M$ can be subdivided into 32 copies of $P^5$. 
The subgroup of symmetries of order 16 of $M$ whose orbit space is $N$ 
acts freely on the 32 copies of $P^5$ subdividing $M$, and so $N$ 
can be described by gluing together two copies of $P^5$. 

Let $\rho$ be the reflection in the second coordinate hyperplane $x_2=0$. 
Then the polytopes $P^5$ and $\rho(P^5)$ are adjacent along a common side. 
We shall describe a side-pairing for the two polytopes $P^5$ and $\rho(P^5)$ 
that glues up to the hyperbolic 5-manifold $N$. 
Table 7 lists the Lorentz normal vectors of the sides of $P^5$ and $\rho(P^5)$.

\begin{table}
$$\begin{array}{lllllll}
s_{1\p{0}}&=&(-1,\p{-}0,\p{-}0,\p{-}0,\p{-}0,\p{-}0) & \quad & s_{17}&=&(-1,\p{-}0,\p{-}0,\p{-}0,\p{-}0,\p{-}0)\\
s_{2\p{0}}&=&(\p{-}0,-1,\p{-}0,\p{-}0,\p{-}0,\p{-}0) & \quad & s_{18}&=&(\p{-}0,\p{-}1,\p{-}0,\p{-}0,\p{-}0,\p{-}0)\\
s_{3\p{0}}&=&(\p{-}0,\p{-}0,-1,\p{-}0,\p{-}0,\p{-}0) & \quad & s_{19}&=&(\p{-}0,\p{-}0,-1,\p{-}0,\p{-}0,\p{-}0)\\
s_{4\p{0}}&=&(\p{-}0,\p{-}0,\p{-}0,-1,\p{-}0,\p{-}0) & \quad & s_{20}&=&(\p{-}0,\p{-}0,\p{-}0,-1,\p{-}0,\p{-}0)\\
s_{5\p{0}}&=&(\p{-}0,\p{-}0,\p{-}0,\p{-}0,-1,\p{-}0) & \quad & s_{21}&=&(\p{-}0,\p{-}0,\p{-}0,\p{-}0,-1,\p{-}0)\\
s_{6\p{0}}&=&(\p{-}0,\p{-}0,\p{-}0,\p{-}1,\p{-}1,\p{-}1) & \quad & s_{22}&=&(\p{-}0,\p{-}0,\p{-}0,\p{-}1,\p{-}1,\p{-}1)\\
s_{7\p{0}}&=&(\p{-}0,\p{-}0,\p{-}1,\p{-}0,\p{-}1,\p{-}1) & \quad & s_{23}&=&(\p{-}0,\p{-}0,\p{-}1,\p{-}0,\p{-}1,\p{-}1)\\
s_{8\p{0}}&=&(\p{-}0,\p{-}0,\p{-}1,\p{-}1,\p{-}0,\p{-}1) & \quad & s_{24}&=&(\p{-}0,\p{-}0,\p{-}1,\p{-}1,\p{-}0,\p{-}1)\\
s_{9\p{0}}&=&(\p{-}0,\p{-}1,\p{-}0,\p{-}0,\p{-}1,\p{-}1) & \quad & s_{25}&=&(\p{-}0,-1,\p{-}0,\p{-}0,\p{-}1,\p{-}1)\\
s_{10}&=&(\p{-}0,\p{-}1,\p{-}0,\p{-}1,\p{-}0,\p{-}1) & \quad & s_{26}&=&(\p{-}0,-1,\p{-}0,\p{-}1,\p{-}0,\p{-}1)\\
s_{11}&=&(\p{-}0,\p{-}1,\p{-}1,\p{-}0,\p{-}0,\p{-}1) & \quad & s_{27}&=&(\p{-}0,-1,\p{-}1,\p{-}0,\p{-}0,\p{-}1)\\
s_{12}&=&(\p{-}1,\p{-}0,\p{-}0,\p{-}0,\p{-}1,\p{-}1) & \quad & s_{28}&=&(\p{-}1,\p{-}0,\p{-}0,\p{-}0,\p{-}1,\p{-}1)\\
s_{13}&=&(\p{-}1,\p{-}0,\p{-}0,\p{-}1,\p{-}0,\p{-}1) & \quad & s_{29}&=&(\p{-}1,\p{-}0,\p{-}0,\p{-}1,\p{-}0,\p{-}1)\\
s_{14}&=&(\p{-}1,\p{-}0,\p{-}1,\p{-}0,\p{-}0,\p{-}1) & \quad & s_{30}&=&(\p{-}1,\p{-}0,\p{-}1,\p{-}0,\p{-}0,\p{-}1)\\
s_{15}&=&(\p{-}1,\p{-}1,\p{-}0,\p{-}0,\p{-}0,\p{-}1) & \quad & s_{31}&=&(\p{-}1,-1,\p{-}0,\p{-}0,\p{-}0,\p{-}1)\\
s_{16}&=&(\p{-}1,\p{-}1,\p{-}1,\p{-}1,\p{-}1,\p{-}2) & \quad & s_{32}&=&(\p{-}1,-1,\p{-}1,\p{-}1,\p{-}1,\p{-}2)
\end{array}$$
\caption{The Lorentz normal vectors of the 32 sides of $P^5$ and $\rho(P^5)$}
\end{table}

We order the sides of $P^5$ and $\rho(P^5)$ by the ordering in Table 7. 
The $i$th side of $\{P^5,\rho(P^5)\}$ is paired to the $j$th side of $\{P^5,\rho(P^5)\}$ 
for $(i,j)=(1,13)$, $(2,18)$, $(3,8)$, $(4,6)$, $(5,20)$, $(7,22)$, $(9,16)$, $(10,19)$,  
$(11,25)$, $(12,29)$, $(14,24)$, $(15,31)$, $(17,30)$, $(21,23)$, $(26,32)$, $(27,28)$. 
The Lorentzian matrices that pair the $i$th side to the $j$th side are given below. 
The second matrix is the identity matrix because side 2 of $P^5$ is 
equal to side 18 of $\rho(P^5)$.

$$\left(\begin{array}{cccccc}
          1 & 0 &          -1 & 0 &\phantom{-}0 & 1 \\
          0 & 0 &\phantom{-}0 & 1 &\phantom{-}0 & 0 \\
          0 & 0 &          -1 & 0 &          -1 & 1 \\
		  1 & 0 &\phantom{-}0 & 0 &          -1 & 1 \\
		  0 & 1 &\phantom{-}0 & 0 &\phantom{-}0 & 0 \\
		  1 & 0 &          -1 & 0 &          -1 & 2 
        \end{array} \right),
\left(\begin{array}{cccccc}
          1 & 0 & 0 & 0 & 0 & 0 \\
          0 & 1 & 0 & 0 & 0 & 0 \\
          0 & 0 & 1 & 0 & 0 & 0 \\
		  0 & 0 & 0 & 1 & 0 & 0 \\
		  0 & 0 & 0 & 0 & 1 & 0 \\
		  0 & 0 & 0 & 0 & 0 & 1 
        \end{array} \right),$$
$$\left(\begin{array}{cccccc}
          0 & 1 & 0 &\phantom{-}0 &\phantom{-}0 & 0 \\
          1 & 0 & 0 &\phantom{-}0 &\phantom{-}0 & 0 \\
          0 & 0 & 1 &          -1 &\phantom{-}0 & 1 \\
		  0 & 0 & 1 &\phantom{-}0 &          -1 & 1 \\
		  0 & 0 & 0 &          -1 &          -1 & 1 \\
		  0 & 0 & 1 &          -1 &          -1 & 2 
        \end{array} \right),
\left(\begin{array}{cccccc}
                    -1 & 0 & 0 & 0 &          -1 & 1 \\
          \phantom{-}0 & 0 & 1 & 0 &\phantom{-}0 & 0 \\
          \phantom{-}0 & 1 & 0 & 0 &\phantom{-}0 & 0 \\
		  \phantom{-}0 & 0 & 0 & 1 &          -1 & 1 \\
		            -1 & 0 & 0 & 1 &\phantom{-}0 & 1 \\
		            -1 & 0 & 0 & 1 &          -1 & 2 
        \end{array} \right),$$
 $$\left(\begin{array}{cccccc}
          0 &          -1 &          -1 &\phantom{-}0 &\phantom{-}0 &\phantom{-}1 \\
          0 &\phantom{-}0 &\phantom{-}1 &\phantom{-}1 &\phantom{-}0 &          -1 \\
          1 &\phantom{-}0 &\phantom{-}0 &\phantom{-}0 &\phantom{-}0 &\phantom{-}0 \\
		  0 &\phantom{-}0 &\phantom{-}0 &\phantom{-}0 &          -1 &\phantom{-}0 \\
		  0 &          -1 &\phantom{-}0 &          -1 &\phantom{-}0 &\phantom{-}1 \\
		  0 &          -1 &          -1 &          -1 &\phantom{-}0 &\phantom{-}2 
        \end{array} \right),
\left(\begin{array}{cccccc}
          0 &    -1 &    -1 &\p{-}0 &\p{-}0 &\p{-}1 \\
          0 &\p{-}0 &\p{-}1 &\p{-}1 &\p{-}0 &    -1 \\
          1 &\p{-}0 &\p{-}0 &\p{-}0 &\p{-}0 &\p{-}0 \\
		  0 &\p{-}0 &    -2 &\p{-}0 &    -1 &\p{-}2 \\
		  0 &    -1 &    -2 &    -1 &    -2 &\p{-}3 \\
		  0 &    -1 &    -3 &    -1 &    -2 &\p{-}4 
        \end{array} \right),$$
$$\left(\begin{array}{cccccc}
              -1 &-2 &\p{-}0 &    -1 &-2 & 3 \\
          \p{-}0 &-2 &\p{-}0 &\p{-}0 &-1 & 2 \\
              -1 &-2 &    -1 &\p{-}0 &-2 & 3 \\
		  \p{-}0 &-1 &\p{-}0 &\p{-}0 &-2 & 2 \\
		  \p{-}0 &-2 &    -1 &    -1 &-2 & 3 \\
		      -1 &-4 &    -1 &    -1 &-4 & 6 
        \end{array} \right),
  \left(\begin{array}{cccccc}
          \p{-}0 &\p{-}0 & 0 &\p{-}0 & 1 &\p{-}0 \\
          \p{-}1 &\p{-}0 & 0 &\p{-}1 & 0 &    -1 \\
          \p{-}0 &\p{-}1 & 0 &\p{-}1 & 0 &    -1 \\
		  \p{-}0 &\p{-}0 & 1 &\p{-}0 & 0 &\p{-}0 \\
		      -1 &    -1 & 0 &\p{-}0 & 0 &\p{-}1 \\
		      -1 &    -1 & 0 &    -1 & 0 &\p{-}2 
        \end{array} \right),$$
 $$\left(\begin{array}{cccccc}
          0 &    -1 &\p{-}0 &    -1 &\p{-}0 &\p{-}1 \\
          0 &\p{-}2 &\p{-}2 &\p{-}1 &\p{-}1 &    -3 \\
          0 &    -1 &\p{-}0 &\p{-}0 &    -1 &\p{-}1 \\
		  1 &\p{-}0 &\p{-}0 &\p{-}0 &\p{-}0 &\p{-}0 \\
		  0 &    -2 &    -1 &\p{-}0 &\p{-}0 &\p{-}2 \\
		  0 &    -3 &    -2 &    -1 &    -1 &\p{-}4 
        \end{array} \right),
  \left(\begin{array}{cccccc}
              -2 &    -1 & 0 &    -1 &    -2 &\p{-}3 \\
          \p{-}0 &\p{-}0 & 0 &\p{-}1 &\p{-}1 &    -1 \\
          \p{-}0 &    -1 & 0 &\p{-}0 &    -1 &\p{-}1 \\
		      -1 &\p{-}0 & 0 &\p{-}0 &    -2 &\p{-}2 \\
		  \p{-}0 &\p{-}0 & 1 &\p{-}0 &\p{-}0 &\p{-}0 \\
		      -2 &    -1 & 0 &    -1 &    -3 &\p{-}4 
        \end{array} \right),$$
 $$\left(\begin{array}{cccccc}
          \p{-}0 &\p{-}0 &\p{-}0 &\p{-}0 & 1 &\p{-}0 \\
          \p{-}1 &\p{-}0 &\p{-}0 &\p{-}1 & 0 &    -1 \\
              -2 &    -1 &    -2 &    -1 & 0 &\p{-}3 \\
		      -2 &\p{-}0 &    -1 &\p{-}0 & 0 &\p{-}2 \\
		      -1 &    -1 &\p{-}0 &\p{-}0 & 0 &\p{-}1 \\
		      -3 &    -1 &    -2 &    -1 & 0 &\p{-}4 
        \end{array} \right),
  \left(\begin{array}{cccccc}
              -1 &    -2 & 0 & 0 & 0 &\p{-}2 \\
          \p{-}2 &\p{-}1 & 0 & 0 & 0 &    -2 \\
          \p{-}0 &\p{-}0 & 1 & 0 & 0 &\p{-}0 \\
		  \p{-}0 &\p{-}0 & 0 & 1 & 0 &\p{-}0 \\
		  \p{-}0 &\p{-}0 & 0 & 0 & 1 &\p{-}0 \\
		      -2 &    -2 & 0 & 0 & 0 &\p{-}3 
        \end{array} \right),$$
  $$\left(\begin{array}{cccccc}
          1 &\p{-}0 &\p{-}0 &    -1 &\p{-}0 & 1 \\
          0 &\p{-}0 &\p{-}0 &\p{-}0 &    -1 & 0 \\
          1 &\p{-}0 &    -1 &\p{-}0 &\p{-}0 & 1 \\
		  0 &    -1 &\p{-}0 &\p{-}0 &\p{-}0 & 0 \\
		  0 &\p{-}0 &    -1 &    -1 &\p{-}0 & 1 \\
		  1 &\p{-}0 &    -1 &    -1 &\p{-}0 & 2 
        \end{array} \right),
  \left(\begin{array}{cccccc}
          \p{-}0 &\p{-}0 &    -1 &    -1 & 0 &\p{-}1 \\
          \p{-}0 &    -1 &\p{-}1 &\p{-}0 & 0 &    -1 \\
              -1 &\p{-}1 &    -1 &    -1 & 1 &\p{-}2 \\
		      -1 &\p{-}0 &    -1 &\p{-}0 & 0 &\p{-}1 \\
		  \p{-}0 &\p{-}0 &    -1 &\p{-}0 & 1 &\p{-}1 \\
		      -1 &\p{-}1 &    -2 &    -1 & 1 &\p{-}3 
        \end{array} \right),$$
  $$\left(\begin{array}{cccccc}
             -1 &\p{-}2 &    -1 &    -2 &\p{-}0 &\p{-}3 \\
         \p{-}0 &    -2 &\p{-}0 &\p{-}1 &\p{-}0 &    -2 \\
         \p{-}0 &\p{-}2 &    -1 &    -2 &    -1 &\p{-}3 \\
		     -1 &\p{-}2 &\p{-}0 &    -2 &    -1 &\p{-}3 \\
		 \p{-}0 &\p{-}1 &\p{-}0 &    -2 &\p{-}0 &\p{-}2 \\
		     -1 &\p{-}4 &    -1 &    -4 &    -1 &\p{-}6 
        \end{array} \right),
  \left(\begin{array}{cccccc}
          \p{-}0 & 1 &    -2 &\p{-}0 &\p{-}0 & 2 \\
              -1 & 0 &\p{-}0 &\p{-}0 &\p{-}0 & 0 \\
          \p{-}0 & 0 &    -1 &    -1 &\p{-}0 & 1 \\
		  \p{-}0 & 0 &    -1 &\p{-}0 &    -1 & 1 \\
		  \p{-}0 & 2 &    -2 &    -1 &    -1 & 3 \\
		  \p{-}0 & 2 &    -3 &    -1 &    -1 & 4
        \end{array} \right).$$

The set $P=P^5\cup\rho(P^5)$ is a right-angled convex polytope with 20 sides. 
The above side-pairing of $\{P^5,\rho(P^5)\}$ determines a facet-pairing 
of $P$ in the sense of \S 11.1 of [7]  
so that $P$ is an inexact fundamental polytope 
for the discrete group $\Gamma$ generated by the above 16 matrices. 
The hyperbolic manifold $N$ is isometric to the space form $H^5/\Gamma$.  

All the entries of the generators of $\Gamma$ are integers, and so 
$\Gamma$ is a torsion-free subgroup of $\Gamma^5$. 
The volume of $N$ is $28\zeta(3)/16=7\zeta(3)/4$. 
The volume of the fundamental domain $\Delta^5$ of $\Gamma^5$ is $7\zeta(3)/15360$, 
and so the index of $\Gamma$ in $\Gamma^5$ is $15360/4=3840$. 
The spherical Coxeter group $(3,3,3,4)$ has order $2^5\cdot 5! = 3840$,  
and so $\Gamma$ is a torsion-free subgroup of $\Gamma^5$ of minimal index. 
Thus $N$ is a minimal volume, integral, hyperbolic 5-manifold. 

R. Kellerhals [3] has proved, by a horoball packing argument,  
that if $M$ is an $m$-cusped hyperbolic 5-manifold, then 
$${\rm vol}(M)/m \ >\  0.3922.$$
Now $N$ has two cusps and 
$$\quad\quad{\rm vol}(N)/2 \ = \ 1.0518\ldots .$$
This suggests that $N$ has very small volume, 
and perhaps $N$ is a minimal volume open hyperbolic 5-manifold. 

\vspace{.2in}
\centerline{References}

\begin{enumerate}
\item B. Everitt, Coxeter groups and hyperbolic manifolds, preprint, 
arXiv:math. GT/0205157 v2, 17 June 2003.

\item J. A. Hillman, Flat 4-manifold groups, 
{\it New Zealand J. Math.} {\bf 24} (1995), 29-40. 

\item R. Kellerhals, Volumes of cusped hyperbolic manifolds, 
{\it Topology} {\bf 37} (1998), 719-734. 

\item R. D. Levine, {\it The Compact Euclidean Space Forms of Dimension Four}, 
Doctoral Dissertation, University of California, Berkeley, 1970. 

\item M. Newman, {\it Integral Matrices}, Pure and Applied Math., vol. {\bf 45}, 
Academic Press, New York, 1972. 

\item L. Potyagailo and E. B. Vinberg, On right-angled Coxeter groups 
in hyperbolic spaces, preprint 2002. 

\item J. G. Ratcliffe, \!{\it Foundations of Hyperbolic Manifolds},
Graduate Texts in Math., vol. {\bf 149}, Springer-Verlag, Berlin, Heidelberg, and
New York, 1994. 

\item J. G. Ratcliffe, Hyperbolic manifolds, In: {\it Handbook of Geometric Topology}, 
Edited by R. J. Daverman and R. B. Sher, North-Holland, 2002, Amsterdam, 899-920. 

\item J. G. Ratcliffe and S. T. Tschantz, 
Volumes of integral congruence hyperbolic manifolds, 
{\it J. Reine Angew. Math.} {\bf 488} (1997), 55-78. 

\item J. G. Ratcliffe and S. T. Tschantz, 
The volume spectrum of hyperbolic 4-manifolds, 
{\it Experimental Math.} {\bf 9} (2000), 101-125. 

\item A. van der Poorten, A proof that Euler missed ... Ap\'ery's proof of 
the irrationality of $\zeta(3)$, {\it Math. Intelligencer} {\bf 1} (1979), 195-203. 

\item E. B. Vinberg, 
Discrete groups generated by reflections in Lobacevskii spaces,
{\it Math. USSR-Sbornik}, {\bf 1} (1967), 429-444. 

\item J. A. Wolf, 
{\it Spaces of Constant Curvature}, 5th edition, Publish or Perish, Wilmington, 1984.  

\end{enumerate}
\end{document}